\pdfoutput=1
\documentclass[11pt]{article}
\usepackage{amssymb}
\usepackage{graphicx}
\usepackage{subcaption}
\usepackage{array}
\parskip \baselineskip
\newcommand{\ds}{\displaystyle}
\newcommand{\R}{\mathbb{R}}

\newcommand{\s}{\subseteq}

\newcommand{\ol}{\overline}

\newcommand{\ra}{\rightarrow}
\newcommand{\Ra}{\Rightarrow}
\newcommand{\es}{\emptyset}
\newcommand{\LRa}{\Leftrightarrow}
\newcommand{\name}[2]{{#1}{\scriptsize{#2}}}

\usepackage[backend=biber,style=vancouver]{biblatex}

\begin{document}

\begin{center}
{\bf Compression with wildcards: The maximum anticliques of  K\"onig-Egarv\'ary graphs}

\name{M}{arcel} \name{W}{ild}
\\Stellenbosch University, Dept. Mathematics,\\ 7600 Stellenbosch, South Africa
\end{center}
\begin{quote}
A{\scriptsize BSTRACT:} {\footnotesize By definition a  K\"onig-Egarv\'ary graph is such that the cardinality of each maximum matching equals the cardinality of each minimum vertex-cover. Each bipartite graph has this property. We unravel the structure
of the family $MACL(G)$ of all maximum anticliques (independent sets) of a  K\"onig-Egarv\'ary graph $G$. This enables us to enumerate $MACL(G)$ not only in output polynomial time, but also in a compressed format (using wildcards). All of this can be viewed as an instance where {\it all} minimum solutions of a submodular function can be enumerated efficiently. }
\end{quote}

{\bf Article type:} Research



\section{Introduction}

Throughout the article, and without further mention, all structures will be finite.
Furthermore "iff" means "if and only if".
Always $G=(V,E)$ is a simple graph with vertex-set $V$ and edge-set\footnote{Formally $E$ is a set of 2-element sets. Thus if $\{a,b\}\in E$ then $\{a,b\}$ is an edge that "connects" vertex $a$ to vertex $b$.} $E$. Vertices $a,b$ are {\it adjacent } if $\{a,b\}\in E$.  An {\it anticlique} $X\s V$ (also known as {\it independent} set or {\it stable} set) is such that no two vertices in $X$ are adjacent.
  Let $Acl(G)$ be the set of all anticliques of $G$, and let $MaxAcl(G)\s Acl(G)$ be the family of all (inclusion-)maximal members of $Acl(G)$.  Since many algorithms exist to generate $MaxAcl(G)$ in output polynomial time, we won't walk that beaten path - but a related one. Specifically, putting 
  $$\alpha(G):=max\{|X|:\ X\in Acl(G)\}$$
  we  strive to capture the structure of
  the family
$$(1)\quad MACL(G):=\{X\in Acl(G):\ |X|=\alpha(G)\}$$
of all {\it maximum} (cardinality) anticliques. Obviously $MACL(G)\s MaxAcl(G)$.
While a new result (Theorem 1) can be obtained for arbitrary graphs $G$, we otherwise focus on $MACL(G)$ for  K\"onig-Egarv\'ary graphs. Serving as tools, Boolean, respectively submodular functions, maximum matchings, and novel wildcards will  feature prominently.

\vspace{0.5cm}
{\bf 1.1} Here comes the section break up. Section 2 briefly widens the horizon by  zooming out from graphs to hypergraphs. For instance what is a hitting-set on the hypergraph level becomes  vertex-cover on the graph level; and the vertex-covers 
are exactly the complements  $V\setminus X$ of the anticliques $X$. Consequently the maximum anticliques (=macls) are the complements of the minimum vertex-covers (=mvcs). Accompanying $MACL(G)$ we write $MVC(G)$ for the family of all mvcs, and will take turns adopting both view points.

As to Section 3, for each graph $G=(V,E)$ we split $V=A^*\uplus F^*\uplus NO^*$ according to whether $v\in V$ belongs to all, some, or no macl of $G$ (in particular $A^*$ is the so called {\it core} of $G$). It turns out that this decomposition is suitable to describe $MACL(G)$ to some extent, even more so when combined with results of  \cite{BML}. 

Section 4 introduces pure 2-CNFs (i.e. Boolean functions in Conjunctive Normal Form all of whose clauses have length two) and indicates how the modelsets $Mod(f)$ of pure 2-CNFs  $f$ can be compressed as disjoint unions of 012-rows like $(0,2,1,1,2)$. Here "2" is the familiar don't-care symbol which in other texts is often  written as $"*"$. Special attention is devoted to the case where all clauses are mixed, i.e. of type $x_i\vee \ol{x}_j$ rather than $x_i\vee x_j$ or $\ol{x}_i\vee\ol{x}_j$. In this scenario $Mod(f)$ can be compressed beyond the use of 012-rows.

Section 5 glimpses at submodular functions $f:2^E\to \R$ because (among other things) they comprise aspects of flows in directed networks $\cal N$. In turn these networks $\cal N$ can be tweeked to render $MVC(B)$ in the special case where $G=B$ is bipartite (due to a result of Picard and Queyranne  \cite{PQ}). It will be insightful to compare this road towards $MVC(B)$ with our quite different approach in Section 8.

 As to Section 6, it is known that the mvcs of a  K\"onig-Egarv\'ary {\it with} a perfect matching and the models of a pure 2-CNF  are essentially the "same thing". We reprove this in a fresh way in Theorem 3. The latter, in tandem with a maximum matching,
   also is the key to unravel (in Section 7) the structure of $MACL(G)$ and $core(G)$ when $G$ is a general K\"onig-Egarv\'ary graph. 
   
   Sections 8 and 9 trim matters in the special cases where the KE graph $G$ is bipartite ($G=B$) respectively a tree ($G=Tr$). In particular, Section 8 shows that, apart from Picard-Queyranne, our approach towards $MVC(B)$ also bears similarities to the Dulmage-Mendelson decomposition of  $B$ into three parts. This decomposition, among other applications (such as disentangling overconstrained systems of equations), classifies all maximum matchings of $B$. In Section 9 we e.g. show that for trees $Tr$ "classifies" improves to an elegant and fast enumeration of all maximum matchings.

\section{Hypergraphs, matchings, and K\"onig-Egarv\'ary graphs}

A {\it hypergraph} on a set $V$ is any family ${\cal H}={\cal H}_V\s {\cal P}(V)$ (=powerset of $V$). The sets $U\in {\cal H}$  are called {\it hyperedges}. Thus a graph is a hypergraph with all hyperedges of cardinality 2. 
In this case we keep on speaking of edges rather than hyperedges and  write $G=(V,E)$ rather than $G=(V,{\cal H})$.
Most graph concepts generalize to hypergraphs.
For instance $X\s V$ is an {\it anticlique} of the hypergraph ${\cal H}_V$ if no hyperedge is contained in $ X$.
Further $Y\s V$ is a {\it hitting-set} of ${\cal H}:={\cal H}_V$ if $Y\cap U\neq \es$ for all hyperedges $U$. 
The  fact that
$$(2)\quad (\forall U\in {\cal H})\ U\not\s X\ \LRa\ (\forall U\in {\cal H})\ (V\setminus X)\cap U\neq \es$$
means that a subset of $V$ is an anticlique of $\cal H$ iff its complement is a hitting-set of $\cal H$. In this correspondence
the {\it maximal} anticliques map upon  the {\it minimal} hitting-sets. It is a notorious open problem whether the family $MHS({\cal H})$ of all minimal hitting-sets (dually: the family of all maximal anticliques) can\footnote{Frankly, more important than formal output-polynomiality is it whether a trillion minimal hitting sets  can be enumerated in a {\it compressed} manner. More on that in 3.4.} be enumerated in output polynomial time \cite{gainerdewar2017}. 

{\bf 2.1} If the hypergraph is a graph $G=(V,E)$ then by (2) the hitting-sets are the {\it vertex-covers (vcs)} of $G$. As opposed to above, in the graph case all minimal vertex-covers (and hence all maximal anticliques)  {\it can} be enumerated in output polynomial time. 
Hence we will turn to greener pastures, i.e. to the enumeration of the family $MVC(G)$ of all minimum-cardinality vertex-covers (=mvcs), and dually to the family $MACL(G)$ (see (1)) of all maximum-cardinality anticliques (=macls). The common cardinality $\tau(G)$ of all $Z\in MVC(G)$ therefore satisfies
$$\tau(G)=|V|-\alpha(G)$$
and hence is generally hard to calculate. For bipartite graphs $G=B$ it evidently holds that $\tau(B)\le min(s_1,s_2)$, where $s_1,s_2$ are the sizes of the {\it shores}\footnote{So $V=S_1\uplus S_2$ and all edges are between $S_1$ and $S_2$. The bipartition $\{S_1,S_2\}$ is unique iff $B$ is connected.} $S_1,S_2$. The inequality can be sharp:

\begin{center}
\includegraphics[scale=1.1]{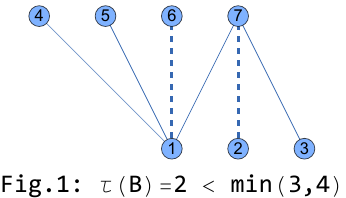}
\end{center}

{\bf 2.2} The hitting-sets of a hypergraph ${\cal H}_V$ are not only related to the anticliques of ${\cal H}_V$ but, equally important, also to the {\it matchings} $M$ of ${\cal H}_V$, i.e. to sets $M\s{\cal H}_V$ of {\it disjoint} hyperedges. It holds that
$$(3)\quad |M|\le |Z| \ for\ all\ matchings\ M\s{\cal H}_V\ and\ all\ hitting\ sets\ Z\s V\ of \ {\cal H}_V.$$
Indeed, if say $M=\{X_1,...,X_t\}$ then $\bigcup M=X_1\uplus\cdots\uplus X_t$, and so 
$$|M|=t\le |Z\cap(X_1\uplus\cdots\uplus X_t)|\le |Z|.$$ 
One says that a hypergraph is {\it K\"onig-Egarv\'ary (KE)} if for some matching $M$ and some hitting-set $Z$ it holds that $|M|=|Z|$.  In this case $M$ has necessarily maximum cardinality among all matchings, and $Z$ has minimum cardinality among all hitting-sets.

{\bf 2.3} So far the KE property has been mainly investigated for graphs, and we will contribute to this endeavour.
Rephrasing (3) for graphs $G=(V,E)$ it reads:
$$(4)\quad |M|\le |Z| \ for\ all\ matchings\ M\s E\ and\ vertex\!-\!covers\ Z\s V\ of \ G$$

\begin{center}
\includegraphics[scale=1.1]{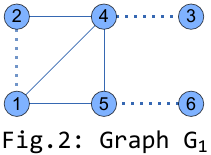}
\end{center}

For instance $G_1$ in Figure 2 has a 3-element matching (dashed edges) and a 3-element vertex-cover $Z_1:=\{1,4,6\}$. Therefore $G_1$ is a KE-graph.
Adopting standard notation let $\nu(G)$ be the maximum cardinality of a matching. Compared to $\tau(G),\alpha(G)$ one can calculate $\nu(G)$ fairly fast. That's because  by  \cite{Micali} a maximum matching of a general graph $G=(V,E)$ can be found in time $O(\sqrt{w}\ m)$, where $w:=|V|$ and $m:=|E|$.
  By $(4)$ it holds\footnote{The gap in $\nu(G)\le\tau(G)$ can be extreme. For instance, in each $2w$-vertex complete graph $G$ one has $\nu(G)=w<2w-1=\tau(G)$.} that 
  $$\nu(G)\le \tau(G),\  and\ \nu(G)=\tau(G)\  takes\ place\ if\!f\ G\ is\ KE.$$
In particular, if $G$ is KE then $\alpha(G)$ and $\tau(G)$ can be calculated in time $O(\sqrt{w}\ m)$. Crucially, each bipartite graph is KE, see  \cite{HV},  \cite{LP}, or  \cite{schrijver2003}. The non-bipartite graph $G_1$ (it has triangles) shows that KE graphs are strictly more general than bipartite graphs.

\section{ $MACL(G)$ and $core(G)$ for general graphs $G$}

Upon discussing the state of the art method to  enumerate $MACL(G)$ for arbitrary graphs $G$, we turn to perfect graphs and reprove that here polynomial total time is achievable. Theorem 1 states that, for any graph $G$, the calculation of   $MACL(G)$ can be reduced to the calculation of $core(G)$ and $MACL(G[F^*])$ for a certain induced subgraph $G[F^*]$. Several results of  \cite{BML} will be woven into our framework.

{\bf 3.1}  Let $G=(V,E)$ be any graph. If $v\in V$ is fixed, then $MACL(G)$ splits into the classes $MACL(G,v)$ and $MACL(G,\ol{v})$ of macls with and without $v$ respectively. Either one, but evidently not both, can be empty. No claim to originality is made for Lemma 1 (yet see 3.3). Adopting standard notation the {\it closed} neihborhood of vertex  $v$ is $NH[v]:=NH(\{v\})\cup\{v\}$.

{\bf Lemma 1:} {\it
\begin{enumerate}
\item[(a)] $MACL(G,\ol{v})\not=\emptyset$ if and only if $\alpha(\ G[V\setminus\{v\}]\ )=\alpha(G)$.
\item[(b)] $MACL(G,v)\not=\emptyset$ if and only if $\alpha(\ G[V\setminus NH[v]]\ )=\alpha(G)-1$.
\end{enumerate} }

{\it Proof.} As to (a), both directions $\Ra$ and $\Leftarrow$ are obvious.

As to $\Ra$ in (b), any $X\in MACL(G,v)\not=\emptyset$ is an anticlique containing $v$, hence $X\cap NH[v]=\{v\}$, hence $X\setminus \{v\}$ is a cardinality $\alpha(G)-1$ anticlique of $G[V\setminus NH[v]\ ]$. It cannot be that $\alpha(G[V\setminus NH[v]])\ge\alpha(G)$ since otherwise $G$ had anticliques $Y\uplus\{v\}$ of cardinality $\ge \alpha(G)+1$, contradiction. Thus $\alpha(G\big[V\setminus NH[v]\big]\ )=\alpha(G)-1$.

 As to $\Leftarrow$ in (b), take any $Y\in MACL(G[V\setminus NH[v]]\ )$, which hence has $|Y|=\alpha(G)-1$.
 Because of $Y\cap NH[v]=\es$ the set $X:=Y\cup\{v\}$ remains an anticlique. From $|X|=|Y|+1=\alpha(G)$ follows
$X\in MACL(G,v)$. Hence $ MACL(G,v)\not =\emptyset$. $\blacksquare$

Let $\cal K$ be a class of graphs  which is closed under taking induced subgraphs.
Suppose there is an algorithm which for each graph $G\in {\cal K}$ with $w$ vertices produces a macl in time $\le T(w)$. Then a recursive search based on Lemma 1 shows that the $N$ macls of any graph $G\in{\cal K}$ can be enumerated in time $O(NwT(w))$. Trouble is, for the class $\cal K$ of {\it all} graphs the state of the art bound for $T(w)$ is $O(1.1996^w)$ \cite{Xiao}. 

If $\cal K$ is the class of all perfect graphs \cite[ch. 65]{schrijver2003} then  $T(w)=O(w^6)$ approximates\footnote{ This is an  AI upper bound which seems to be close to the current (implicite) record; apparently none of several approaches has yet claimed to be the record holder.} the current state of the art (and so the macls of perfect graphs can be calculated in output polynomial time). For many subclasses (such as bipartite, chordal, and comparability graphs) the $O(w^6)$ bound can be trimmed considerably.

 {\bf 3.2} Recall that $core(G)$ is the intersection of all macls of $G=(V,E)$. 
 We usually put $A^*(G):=core(G)$ since this insinuates\footnote{The superscript of $A^*$ will make sense in 3.4.} that $A^*=A^*(G)$ is the set of vertices appearing in {\it all} macls. Likewise $NO^*$  is the set of vertices appearing in {\it no} macl. Further $F^*:=V\setminus (A^*\cup NO^*)$ is the set of all {\it flexible} vertices $v$ in that $v$ sits in some but not all macl's. Thus $V=A^*\uplus F^*\uplus NO^*$.

\begin{center}
\includegraphics[scale=1.1]{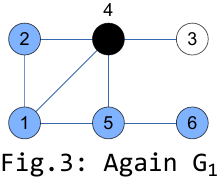}
\end{center}

For instance, one checks that  $G_1=(V_1,E_1)$ from Figure 2 has\\ $MACL(G_1)=\big\{\{1,3,6\},\{2,3,5\},\{2,3,6\}\big\}$, and so $V_1=A^*\uplus F^*\uplus NO^*$ is visualized by the white, blue, black vertices in Figure 3.

  For each $U\s V$ the {\it neighborhood} $NH(U)$ consists of all vertices which are adjacent to some vertex in $U$. Put another way $NH(U)$ is the union of all sets $NH(u)$ when $u$ ranges over $U$. We claim that

$$(5)\quad NH(A^*)\s NO^*,\ hence\  NH(F^*)\s F^*\uplus NO^*\ and\ NH(NO^*)\s A^*\uplus F^*\uplus NO^*.$$

The second and third inclusion in (5) are immediate consequences of the first. As to the first inclusion, $NH(A^*)\cap A^*=\es$ since otherwise there was an edge between members of $A^*$, which is impossible within an anticlique. Similarly, suppose there was an edge $\{u,v\}$ of $G$ with $u\in A^*$ and $v\in F^*$. By definition of $F^*$ there is a macl $X$ of $G$ with $v\in X$. But also $u\in X$ because $u$ is in {\it all} macls. This leads to the contradiction that the anticlique $X$ contains an edge $\{u,v\}$. This shows that besides
$NH(A^*)\cap A^*=\es$ also $NH(A^*)\cap F^*=\es$ takes place, and so $NH(A^*)\s NO^*$ in (5).

{\bf Theorem 1 (the $A^*-F^*-NO^*$ theorem): }{\it For each graph $G$ it holds that 
$MACL(G)\stackrel{!}{=}\ {\cal S}:=\{Y\uplus A^*:\ Y\in MACL(G[F^*])\ \}$.}

{\it Proof.} Take any member $Y\uplus A^*$ of ${\cal S}$. It follows from $NH (A^*)\cap F^*=\es$ (see (5)) that $Y\uplus A^*$ is an anticlique in the first place. From $|Y\uplus A^*|=\alpha(G[F^*])+|A^*|$ follows $\alpha(G)\ge \alpha(G[F^*])+|A^*|$. Conversely,  pick $X\in MACL(G)$. From $V=A^*\uplus F^*\uplus NO^*$ and $X\cap NO^*=\es$ follows that $X=Y'\uplus A^*$ for some $Y'\in Acl(G[F^*])$. Hence $\alpha(G)=|X|=|Y'|+|A^*|\le \alpha(G[F^*])+|A^*|$. One concludes that
$\alpha(G)=\alpha(G[F^*])+|A^*|$ and that $MACL(G)={\cal S}$.
 $\square$

The $A^*-F^*-NO^*$ theorem becomes the more relevant  the easier $A^*$ can be calculated and the more inviting the structure of $G[F^*]$ is. That comes to the fore in Section 8.

{\bf 3.3} Let us point out related results of \cite{BML}. To begin with, up to terminology the first inclusion of (5) was independently established in \cite[Lemma 2]{BML}. In tandem with $core(G)$ the authors define $corona(G)$ as the set of vertices that belong to {\it some} macl; in other words $corona(G)=V\setminus NO^*$. While I could not find an exact reference for Lemma 1 above, similar observations were made in \cite{BML}:
\begin{itemize}
    \item[(6)] $v\in core(G)\LRa \alpha(G[V\setminus\{v\}])=\alpha(G)-1$  \cite[Lemma 4]{BML}
    \item[(7)] $v\in corona(G)\LRa \alpha(G[V\setminus NH[v]\ ])=\alpha(G)-1$  \cite[Lemma 6]{BML}
\end{itemize}

From $\alpha(G)+\tau(G)=|V|$ and $\nu(G)\le \tau(G)$ (Section 2) follows $\alpha(G)+\nu(G)\le|V|$, but generally neither $\alpha(G)\le \nu(G)$ nor $\alpha(G)\ge \nu(G)$ takes place. If however $G$ is connected and $\alpha(G)> \nu(G)$ then 
$$|core(G)|\ge 1+\alpha(G)-\nu(G).$$
This inequality is sharp \cite[Thm.1]{BML} and improves previous lower bounds of $|core(G)|$. Here comes an interesting property of the subgraph $G[F^*]$ appearing in Thm.1:
\begin{itemize}
\item[(8)] {\it Each $Y\in MACL(G[F^*])$ can be matched to a disjoint $Y_0\in MACL(G[F^*])$. }
\end{itemize}

To spell it out, for each $Y\in MACL(G[F^*])$ there is a matching $M_0=\{\{v_1,v_1'\},...,\{v_t,v_t'\}\}$ of $G[F^*]$ such that $Y=\{v_1,...,v_t\}$ and $Y_0=\{v'_1,...,v'_t\}$.
For instance, one checks that the pentagon has five macls $Y$, and each one has a  such a companion $Y_0$.

{\it Proof of (8).} By Theorem 1 each  $Y\in MACL(G[F^*])$ extends to the macl $S:=Y\uplus A^*$ of $G$. By [BML,Lemma 3] there is a matching $M_0$ from $S\setminus core(G)$ into $corona(G)\setminus S$. But $S\setminus core(G)$ is $Y$, and $Y_0:=corona(G)\setminus S$ is contained in $G[F^*]$. $\square$

According to \cite[Thm.2]{BML} the problems of computing the numbers $\alpha(G)$ and $|core(G)|$ are polynomially equivalent. This says nothing about the complexities of {\it listing} $MACL(G)$ respectively $core(G)$. However, 
for perfect graphs $G$ one {\it can} list $MACL(G)$  in output polynomial time (3.1); and by other reasons  \cite[Remark 1]{BML} also $core(G)$ {\it can} be calculated in polynomial time.

{\bf 3.4}    For each graph $G=(V,E)$ it holds that $V=A\uplus F\uplus NO$, where by definition $A$ is the set of vertices appearing in {\it all} mvcs, $NO$ is the set of vertices appearing in {\it no} mvc, and $F:=V\setminus (A\uplus NO)$ is the set of all {\it flexible} vertices $v$ in that $v$ sits in some but not all mvcs. Furthermore
$$A^*=\{v\in V:\ (\forall X\in MACL(G))(v\in X)\}=\{v\in V:\ (\forall X\in MACL(G))(v\not\in X^c)\}$$
$$=\{v\in V:\ (\forall Y\in MVC(G))(v\not\in Y)\}=NO.$$
Similarly one shows that $NO^*=A$. Therefore $F^*=F$ and hence (see also Fig.4) Theorem 1 dualizes as follows:

{\bf Theorem 1' (the $A-F-NO$ theorem): }{\it For each graph $G$ it holds that 
$MVC(G)=\{Z\uplus A:\ Z\in MVC(G[F])\ \}$.}

\begin{figure}[h]
\centering
\includegraphics[scale=0.4]{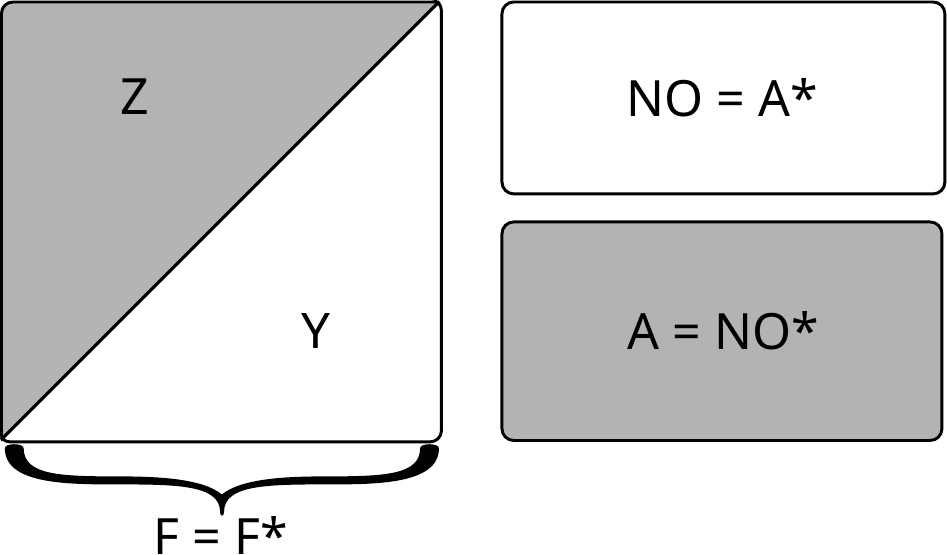}
\caption*{Figure 4: The duality between Theorem 1 and Theorem 1'}
\end{figure}

For a hypergraph $\cal H$ let $MCHS({\cal H})\s MHS({\cal H})$ be the subfamily of all minimum-cardinality hitting sets. While \cite{Wild2} offers no structural analysis of $MCHS({\cal H})$, computational experiments show that $MCHS({\cal H})$ can be compressed even more efficiently (and in more clear-cut ways) than $MHS({\cal H})$. If $\cal H$ is a graph $G$ then $MCHS({\cal H})=MVC(G)$. The article in front of you offers both a structural analysis of $MVC(G)$ (that goes beyond Theorem 1' as $G$ specializes to KE, bipartite, or forest), as well as a compressed enumeration (which however differs a lot from the one for $MCHS({\cal H})$).

\section{A first look at pure Boolean 2-CNFs}

We assume a basic familiarity with Boolean formulas and functions, e.g. see \cite{crama2011}.
The {\it modelset} of a Boolean function $f:\{0,1\}^n\to\{0,1\}$ is $Mod(f):=\{y\in\{0,1\}^n:\ f(y)=1\}$.
One calls $f$ {\it satisfiable} if $Mod(f)\neq\es$. We will express parts of $Mod(f)$ by {\it 012-rows} such as
$$r:=(0,{\bf 2},0,1,{\bf 2},1,1):=\big\{(0,{\bf 0},0,1,{\bf 0},1,1),(0,{\bf 0},0,1,{\bf 1},1,1),(0,{\bf 1},0,1,{\bf 0},1,1),(0,{\bf 1},0,1,{\bf 1},1,1)\big\}.$$
Hence "2" is a dont-care symbol which in other texts is often written as $"*"$. Generally a 012-row with $s$ don't-cares has cardinality $2^s$.
Two Boolean formulas $f,g$ may define Boolean functions with the same modelset. In this case we write $f\equiv g$ and say that $f$ and $g$ are {\it equivalent}.
So each Boolean {\it formula} yields a unique Boolean {\it function} (and that's our way to define them), but not conversely.
Recall that a {\it clause} is a disjunction of Boolean literals, such as $\ol{x}_1\vee x_3\vee \ol{x}_7$. A Boolean formula is in {\it Conjunctive Normal Form (CNF)} if it is a conjunction of clauses. 

{\bf 4.1} If all clauses of a CNF $f:\{0,1\}^n\to\{0,1\}$ have exactly\footnote{In an ordinary 2-CNF also clauses are allowed that consist of a single literal. Such 2-CNFs are easily reduced to pure 2-CNFs, see e.g. \cite[Sec.1 4]{wild2025}.}  two literals then on speaks of a {\it pure 2-CNF}.
The clauses are either {\it positive} ($x_i\vee x_j$) or {\it negative} ($\ol{x}_i\vee \ol{x}_j$) or {\it mixed} ($x_i\vee \ol{x}_j$). Provided $f$ is satisfiable (which can be decided in $O(n+m)$ time where $m$ is the number of 2-clauses), a somewhat technical\footnote{For instance, the first step is to apply so called Horn-renaming in order to get rid of all positive clauses.} algorithm is applicable, which usually (no guarantee) renders $Mod(f)$ in a {\it compressed format}, i.e. as disjoint union of $R$ many 012-rows, and which runs in output  polynomial time $O(Rn^3)$.
This algorithm will be put to good use in Section 6.

{\bf 4.2} If at the outset $f$ has only mixed clauses (this happens in Section 8) then mentioned technicalities disappear and moreover additional wildcards (apart from the don't-care 2's) can be introduced. Let us give the gist of it. Consider this conjunction of mixed 2-clauses:
 
$$(9)\ f_3':=(\ol{y}_1\vee y_4)\wedge(\ol{y}_1\vee y_5)\wedge(\ol{y}_1\vee y_6)\wedge(\ol{y}_6\vee y_7)\wedge(\ol{y}_6\vee y_9)\wedge(\ol{y}_9\vee y_{10})\wedge(\ol{y}_9\vee y_{11}).$$
Recalling the distributive law and the fact that (say) $\ol{y}_1\vee y_4\equiv y_1\ra y_4$, one concludes that
$$(10)\quad  f_3'\equiv\big(y_1\ra(y_4\wedge y_5\wedge y_6)\big)\ \wedge\ \big(y_6\ra(y_7\wedge y_9)\big) \ \wedge\ \big(y_9\ra(y_{10}\wedge y_{11})\big).$$
Following \cite{wild2025} the {\it (a,b)-wildcard} is defined as $(a,b,b,...,b):=({\bf 0},2,2,...,2)\uplus ({\bf 1},1,1,...,1)$. In particular $|(a,b,b,...,b)|=2^s+1$, where $s$ is the number of $b$-symbols.
The {\it (a,b)-algorithm} of \cite{wild2025} renders the modelset of $f_3'$ in (10) as the disjoint union of the {\it 012ab-rows} $\sigma_1,\sigma_2,\sigma_3$ in Table 1:

 \begin{tabular}{c|c|c|c|c|c|c|c|c|l} \hline
 & $y_1$ & $y_4$& $y_5$ & $y_6$ & $y_7$ & $y_9$ &$y_{10}$ &$y_{11}$& cardinality\\ \hline
 &   &  &   &   &  & & & \\ \hline
 $\sigma_1=$ & $a$ & $b$ & $b$ & $1$ & $1$ & $1$ & $1$ & $1$ &   $2^2+1$\\ \hline
 $\sigma_2=$ & $0$ & $2$ & $2$ & $0$ & $2$ & $0$ & $0$ & $2$ &   $2^4$\\ \hline
  $\sigma_3=$ & $0$ & $2$ & $2$ & $0$ & $2$ & $a$ & $1$ & $b$ &   $2^3\cdot (2^1+1)$\\ \hline 
\end{tabular}

{\sl Table 1: Expressing $Mod(f_3')$ in terms of 012ab-rows}

For instance, this is why $\sigma_3\s Mod(f_3')$. Each bitstring $y=(y_1,y_4,...,y_{11})$ in $\sigma_3$ {\it trivially satisfies} $y_1\ra(y_4\wedge y_5\wedge y_6)$ and $y_6\ra(y_7\wedge y_9)$ in the sense that $y_1=0$ and $y_6=0$ respectively. Furthermore $y_{10}=1$   reduces  $y_9\ra(y_{10}\wedge y_{11})$ to $y_9\ra y_{11}$. The latter holds because of the wildcard $(a,b)$ in $\sigma_3$. Each 012ab-row with $s$ many\footnote{They must be distinguished by subscripts as in $(1,a_1,b_1,b_1,0,b_2, 1,2,a_3,a_2,b_3)$.} $(a,b)$-wildcards  can be, if necessary (the prize being a loss of compression), rewritten as disjoint union of $2^s$ many 012-rows.

{\bf 4.3} The $(a,b)$-algorithm shines whenever a family of implications $ x_i\ra(x_j\wedge...\wedge x_k)$ with {\it singleton}\footnote{ General {\it implications} look like  $ (x_i\wedge...\wedge x_\ell) \ra(x_j\wedge...\wedge x_k)$.} premises  needs to be processed. In particular the (a,b)-algorithm is fit to represent the set $Id(P,\le)$ of all order ideals of a poset $(P,\le)$ as a disjoint union of 012ab-rows (which is highly relevant in 5.2). Here each non-minimal $\alpha\in P$, say with lower covers $\beta,\gamma,\delta$, contributes the implication $\alpha\ra\beta\wedge\gamma\wedge\delta$.

\section{Minimizing submodular functions}

In 5.3 and 5.4 we show that in principle our quest (in Sec.8) to compress $MVC(B)$ for bipartite graphs $B$ can be reduced to an application of a theorem in \cite{PQ}. Yet "in principle" is not the same (5.5) as "in practice". Actually we commence with the even more abstract framework of distributive lattices (5.1) and submodular functions (5.2).

{\bf 5.1} The reader is assumed to be familiar with basic concepts concerning posets and lattices. Fact (11) is copied almost literally from \cite[p.128]{Fu2008}:

\begin{itemize}
\item[(11)] {\it Let $\cal D$ be a set of subsets of a finite set $E$ with $\es,E\in{\cal D}$.
Then $\cal D$ is a distributive set-lattice (i.e. the lattice operations are set union and set intersection) iff there exists a poset $(\Pi(E),\le)$ on a partition $\Pi(E)$ of $E$ such that $\cal D$ is expressed as follows:
$$X\in{\cal D}\  \LRa\   there\ exists\ an\ ideal\ {\cal J}\ of\ (\Pi(E),\le) \ such\ that\ X =\bigcup\{F:\ F\in {\cal J}\}.$$ }
\end{itemize}

Following \cite[p.129]{Fu2008} we denote the poset $(\Pi(E),\le)$ appearing in (11) by $\cal P(D)$. Conversely
for any poset $\cal P$ on a partition of $E$ there exists by (11) a unique distributive lattice ${\cal D}\s 2^E$ such that $\cal P=P(D)$. We denote this lattice by $\cal D(P)$.
To summarize,  (11)  characterizes distributive set-lattices ${\cal D}\s 2^E$ which are {\it normalized} in the sense that $\es$ is the smallest and $E$ the largest element of $\cal D$.

{\bf 5.2} Let us turn to distributive set-lattices ${\cal D'}\s 2^{E'}$ which are not necessarily normalized, and which support a function $f:{\cal D'}\to\R$ which is {\it submodular}, i.e.
$$(12)\quad  f(X\cap Y)+f(X\cup Y)\le f(X)+f(Y)\ for\ all\ X,Y\in {\cal D'}.$$

Many combinatorial optimization problems can be reduced to the minimization of  submodular functions; see \cite{Fu2008} and the references therein. Put  $\mu:=min\{f(Y):\ Y\in {\cal D'}\}$ and define the {\it minimizer} of $f$ as 
$${\cal D'}_{min}(f):=\{X\in {\cal D'}:\ f(X)=\mu\}.$$
Knowing the whole of ${\cal D'}_{min}(f)$ is often much better than knowing only one
$X\in {\cal D'}_{min}(f)$. For instance, scanning ${\cal D'}_{min}(f)$ allows to optimize relevant secondary properties of an otherwise "random"  minimum solution.
 The following observation (apparently by Ore 1956) is crucial:
$$(13)\quad {\cal D'}_{min}(f)\ is\ a\ sublattice\ of\ {\cal D'}.$$
Indeed, if $X,Y\in {\cal D'}_{min}(f)$ then $\mu\le f(X\cap Y),f(X\cup Y)$, hence $2\mu\le f(X\cap Y)+f(X\cup Y) \le f(X)+f(Y)=2\mu$, hence $f(X\cap Y)=f(X\cup Y)=\mu$, hence $X\cap Y,X\cup Y\in {\cal D'}_{min}(f).\ \square$

Being a sublattice of a distributive lattice,  ${\cal D'}_{min}(f)$ is distributive as well. Let $E_{min},E_{max}$ be the smallest and largest element of
${\cal D'}_{min}(f)$ respectively. Putting 
$$E:=E_{max}\setminus E_{min}\ and\ {\cal D}:=\{X\setminus E_{min}:\ X\in {\cal D'}_{min}(f)\},$$
it follows that ${\cal D}\s 2^E$ is a normalized distributive set-lattice.  Hence
$$(14)\quad {\cal P}_{min}(f):={\cal P(D)}$$
is well-defined by 5.1 and relates to $f$ in a crucial way: ${\cal P}_{min}(f) $ is a poset of type $(\Pi(E),\le)$ such that the ideals of $(\Pi(E),\le)$, padded  by $E_{min}$, are exactly the members of ${\cal D'}_{min}(f)$. That's why we call ${\cal P}_{min}(f)$ the {\it principal structure}\footnote{This name was used by Fushige in the 80s but he later changed it to "principal partition" (e.g. in \cite{Fu2008}). This is a bit unfortunate since the partial ordering of $\Pi(E)$ is just as important as its blocks.} of the submodular function $f:{\cal D'}\to\R$.

{\bf 5.3} Let us turn to a type of submodular function 
which relates (albeit in hidden ways) to the enumeration of $MVC(B)$ for bipartite graphs $B$. For this the reader is assumed to be familiar with flows in directed networks and the maxflow-mincut theorem, see e.g. \cite{schrijver2003} or\footnote{This reference is recommended for readers who love extended toy examples, be they tailor-made or from real life. Many of them are accompanied with detailed Mathematica code.} \cite{HV}.
Thus let ${\cal N}=(V,A)$ 
be a directed $s-t$ network with vertex-set $V$, arc-set $A$, source $s$, and sink $t$. Recall that a subset $U\s V$ that contains $s$ is called a {\it cut}, and that $\delta^+(U)$ is defined as the sum of the capacities of the arcs with tail in $U$ and head in $V\setminus U$. Consider the (not normalized) distributive set-lattice ${\cal D'}:=\{U\in 2^V:\ s\in U\}$. It is well-known that
$$(15)\quad \delta^+: {\cal D'}\to\R\ \ is\ submodular.$$
Here comes the short proof. The number (say) $a$ in Figure 5(a) is the sum of all capacities of arcs that go from $X$ to $V\setminus (X\cup Y)$. Similarly  $e$ is the sum of all capacities of arcs that go from $X\cap Y$ to $Y\setminus (X\cap Y)$, and so forth. It follows that
$$\delta^+(X\cap Y)+\delta^+(X\cup Y)=(c+d+e)+(a+b+c)$$
$$\le (a+c+e+g)+(b+c+d+h)=\delta^+(X)+\delta^+(Y). \square$$
The principal structure of $\delta^+$ was found in \cite{PQ}. For instance, 
the principal structure of the submodular function $\delta^+$ coupled to the network in Fig.5(b) (taken from \cite{PY}) is rendered in Fig.5(c).

\begin{figure}[htbp]
    \centering
    \begin{subfigure}[b]{0.24\textwidth}
        \centering
        \includegraphics[width=\linewidth]{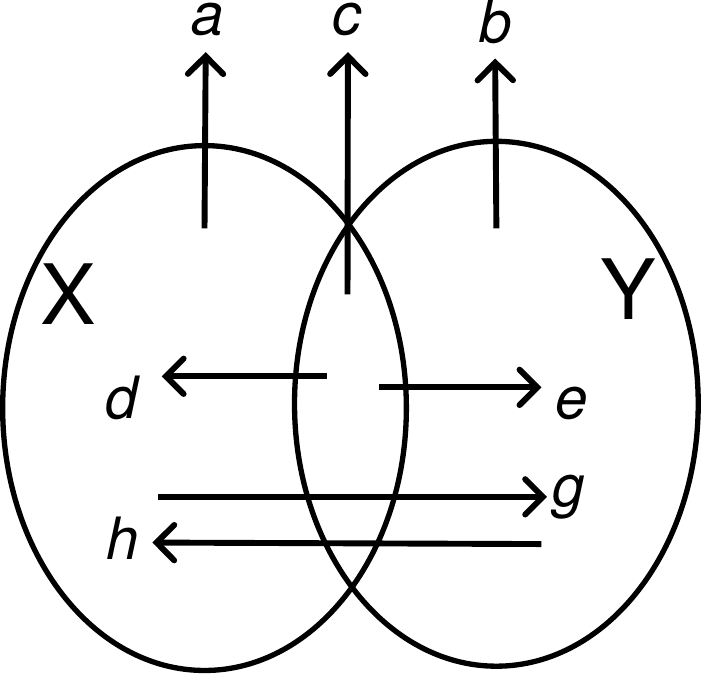}
        \caption*{5(a)}
    \end{subfigure}
    \hfill 
    \begin{subfigure}[b]{0.56\textwidth}
        \centering
        \includegraphics[width=\linewidth]{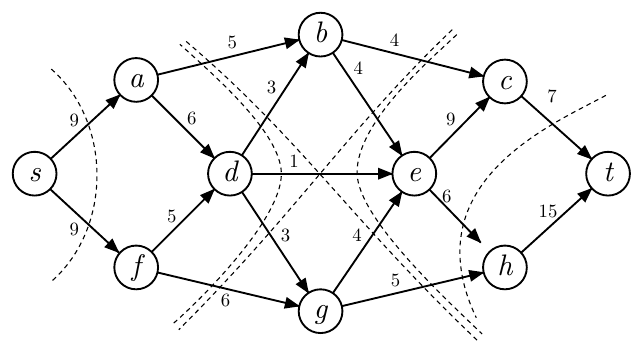}
        \caption*{5(b)}
    \end{subfigure}
    \hfill 
    \begin{subfigure}[b]{0.18\textwidth}
        \centering
        \includegraphics[width=\linewidth]{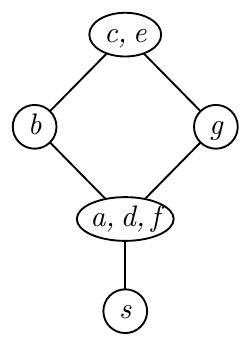}
        \caption*{5(c)}
    \end{subfigure}
   
    \caption*{Figure 5: A directed $s-t$ network and its principal structure.}
\end{figure}

{\bf 5.4}  Given any bipartite graph $B$  let us associate with it an exquisitely chosen network that incorporates the whole of $MVC(B)$.

{\bf Theorem 2:} {\it There is a  bijective correspondence between the  mvcs of $B=(V,E)$ and the mincuts of a suitable $s-t$ network ${\cal N}_B$. }

{\it Proof.} The network ${\cal N}_B$ has vertex-set $\{s,t\}\uplus S\uplus T$, where $S,T$ are the shores of $B$ and where $s,t$ are "artificial" source and sink respectively. By definition the arcset $A$ has three types of arcs. The type (i) arcs go from $s$ to all vertices in $S$, the type (iii) arcs go from the vertices in $T$ to $t$, and the type (ii) arcs are the edges of $E$, except they are directed from $S$ to $T$. The capacity function $c:A\to \R$ assigns value 1 to all type (i) or (iii) arcs, and a large value $L$ (larger than $|V|$) to all type (ii) arcs.

Let $Z:=S_2\uplus T_1$ be an arbitrary mvc of $B$, see Fig.6(a). Putting $S_1:=S\setminus S_2$ and $T_2:=T\setminus T_1$, we claim that $U:=\{s\}\uplus S_1\uplus T_1$ is a mincut of ${\cal N}_B$. First, no arcs starting in $S_1$ are leaving $U$ because, $S_2\uplus T_1$ being a vertex-cover of $B$, the three boldface arrows in Fig.6(a) indicate  the only type (ii) arcs present in 
${\cal N}_B$.   Consequently
$\delta^+(U)=|S_2|+|T_1|$.
Next, since $B$ is KE, there is a matching $M=\{\{v_1,v_1'\},\{v_2,v_2'\},...,\{v_\tau,v_\tau'\}\}$ of $B$ with all $v_i\in S$, all $v_i'\in T$, and with $|M|=\tau=|Z|$. 
This matching triggers a flow $f\ell$ of ${\cal N}_B$ as follows.
For all type (ii) arcs put $f\ell(v,v'):=1$ if $\{v,v'\}\in M$, and $f\ell(v,v'):=0$
otherwise. Further, for all $v\in S$ put $f\ell(s,v):=1$ if $v=v_i$ for some index $i$, and $f\ell(s,v):=0$ otherwise. Likewise  for all $v'\in T$ put $f\ell(v',t):=1$ if $v'=v_i'$ for some index $i$, and $f\ell(v',t):=0$ otherwise.
By definition the value of a flow is the sum of the values of the arcs with tail $s$.
 Hence $val(f\ell)=|\{v_1,...,v_\tau\}|=|Z|=\delta^+(U)$. Therefore (by the trivial direction in the Ford-Fulkerson Theorem) the  cut $U$ must be minimum.

\begin{figure}[h]
\centering
\includegraphics[scale=0.5]{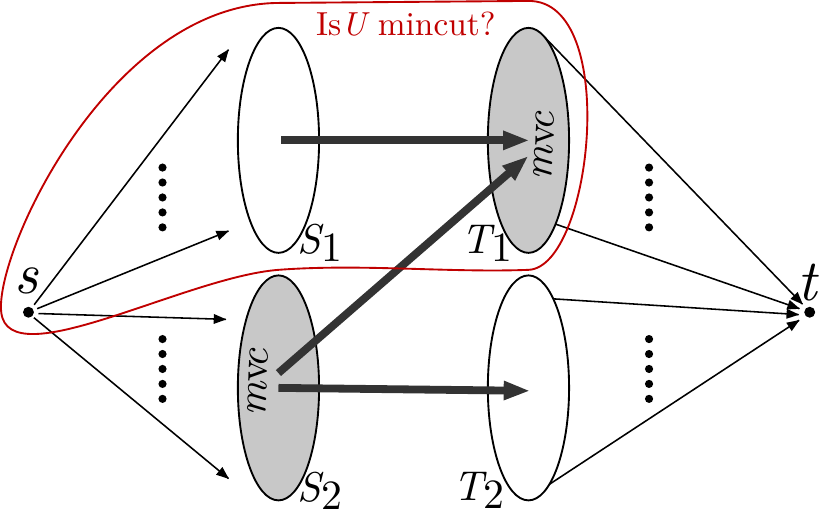}
\caption*{Fig.6(a): Each mvc yields a mincut}
\end{figure}

Conversely,  let $\ol{U}$ be a minimum cut of ${\cal N}_B$. In order to pinpoint some mvc $Z'$ that induces $\ol{U}$, put
$S_1':=\ol{U}\cap S,\ T_1':=\ol{U}\cap T,\ S_2':=S\setminus S_1',\ T_2':=T\setminus T_1'$; see Figure 6(b). Obviously $\delta^+(\ol{U})\ge |S_2|$ but our aim is to show
$\delta^+(\ol{U})= |S_2|+|T_1'|$.
No type (ii) arcs go
 from $S_1'$ to $T_2'$ since otherwise $\delta^+(\ol{U})\ge L$, which contradicts $\ol{U}$ being a minimum cut. Hence all arcs with tail in $S_1'$ have the head in $T_1'$. Hence in $B$ we have $NH(S_1')\s T_1'$. If there was a  vertex $v'\in T_1'\setminus NH(S_1')$ then the nature\footnote{There may well be arcs of type $(v,v')$ but then by assumption $v\in S_2'$. Therefore  $(v,v')$ is not an arc with tail in 
  $\ol{U}\setminus\{v'\}$, and hence cannot increase the value of $\delta^+(\ol{U}\setminus\{v'\})$ beyond $\delta^+(\ol{U})-1$ .} of $v'$ implies that $\delta^+(\ol{U}\setminus\{v'\})=\delta^+(\ol{U})-c(v',t)=\delta^+(\ol{U})-1$. 
 This again contradicts the minimality of $\ol{U}$. Therefore $NH(S_1')= T_1'$. This implies that $\delta^+(\ol{U})=|S_2'|+|T_1'|$, and that $Z':=S_2'\uplus T_1'$ is a vertex-cover of $B$. But is $Z'$ of minimum cardinality? By Ford-Fulkerson there is a flow $\ol{f\ell}$ with $val(\ol{f\ell})=\delta^+(\ol{U})=|Z'|$. Since all capacities of type (i) or (iii) arcs are 1, since $\ol{f\ell}$ is integer-valued, and since the conservation laws hold, the type (ii) arcs with $\ol{f\ell}(v,v')=1$ must be such that they constitute (upon undirecting them) a matching $\ol{M}$ of $B$.
 Therefore $|\ol{M}|=val(\ol{f\ell})$, i.e. $|\ol{M}|=|Z'|$. The trivial part of the K\"onig Egarv\'ary Theorem implies that $Z'$ is a mvc.  $\square$

\begin{figure}[h]
\centering
\includegraphics[scale=0.5]{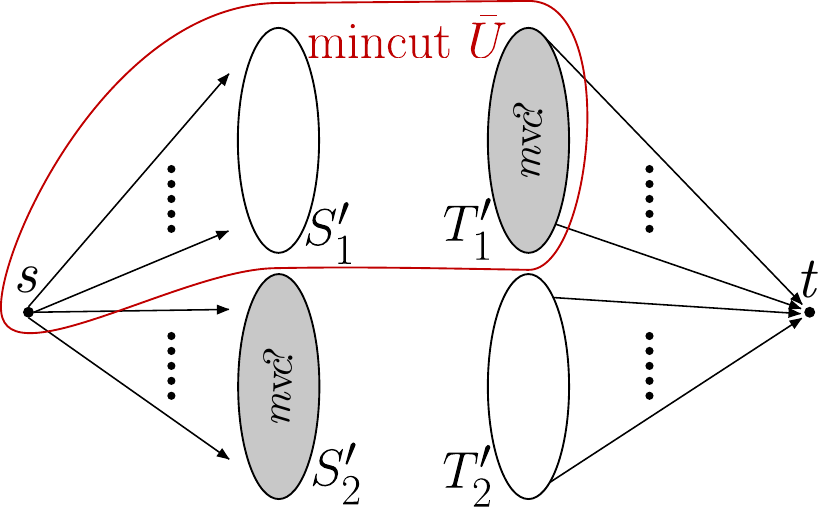}
\caption*{Fig.6(b): Each mincut yields a mvc}
\end{figure}

{\bf 5.5} We see that the task to find the family $MVC(B)$ of all mvcs of a bipartite graph $B$ can in principle be solved by applying the algorithm in \cite{PQ} to a very special type of $s-t$ network. However, using quite different ideas (e.g. Boolean functions and wildcards) we will calculate $MVC(B)$ faster, and in compressed format.

As a preview of things to come consider the bipartite graph $B_2$ in Fig 7a. Let $\delta_2^+$ be the submodular function that is coupled (as explained in 5.4) to the network ${\cal N}_{B_2}$. It then holds (details omitted for the time being) that the principal structure ${\cal P}_{min}(\delta_2^+)$ is rendered by Fig.7(b). Hence
$MVC(B_2)$ is bijective\footnote{One could thus enforce a distributive lattice strucure directly on $MVC(B_2)$. However the operations in this lattice would be more intricate than just set union and set intersection.} to the ideals of the poset in Fig.7(b).   For instance, the ideal $I_1=\{5,6,7\}$ yields the  mvc $\{\ol{1},\ol{2},\ol{3},\ol{4},5,6,7\}$. Similarly $I_2=\{1,2,5\}$ yields the 
mvc $\{1,2,5,\ol{3},\ol{4},\ol{6},\ol{7}\}$. The union $I_1\cup I_2$ yields the
mvc $\{1,2,5,6,7,\ol{3},\ol{4}\}$.  Altogether there are 26 ideals and whence 26 mvcs. As always, the principal structure is isomorphic to the subposet 
$\big\{\{1\}\!\downarrow,\{2\}\!\downarrow,...,\{7\}\!\downarrow\big\}$ of join-irreducibles of its ideal-lattice. Recall that this ideal lattice is the minimizer ${\cal D}'_{min}(\delta_2^+)$, see Fig.7(c).

\begin{center}
\includegraphics[scale=0.95]{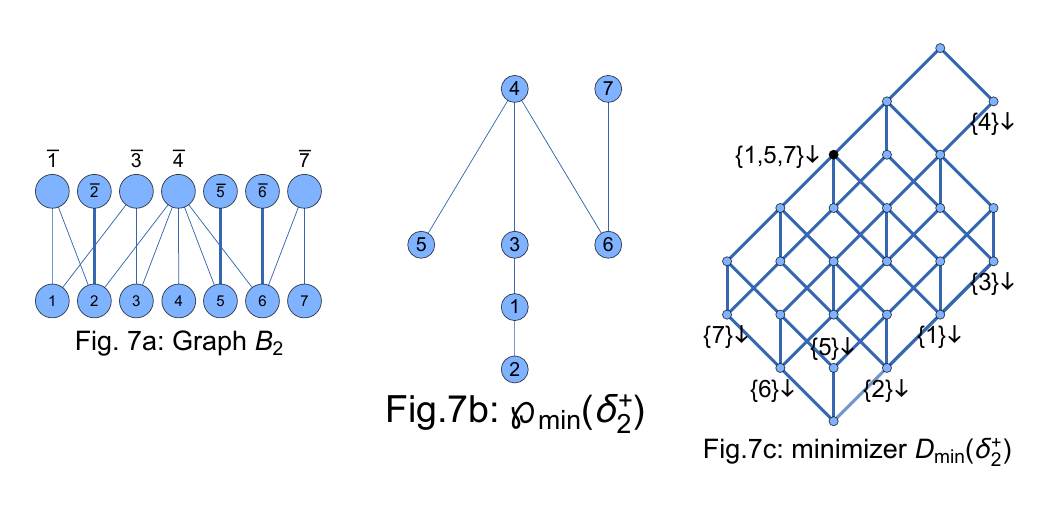}
\end{center}

\section{Compressing $MVC(G)$ for KE graphs $G$\\ with a perfect matching}

 K\"onig-Egarv\'ary graphs with perfect matching are less haphazard than they may  seem. They pave the way to understanding the structure of $MACL(B)$ for general bipartite graphs $B$.
  Recall from Section 3 the duality between $MACL(G)$ and $MVC(G)$. In Section 6 it suits us better to stick to $MVC(G)$.
With each pure and satisfiable 2-CNF one can associate a certain "matched graph". Such graphs  are always K\"onig-Egarv\'ary  (KE) and furthermore possess a perfect matching. The converse is true as well.
 For instance $G_1$ in Fig.2 is KE with perfect matching  and hence is the matched graph of a certain pure and satisfiable 2-CNF.
 The following device which was used in previous publications will come in handy. Loosely speaking the $g$-{\it wildcard} $(g,g,...,g)$ means "exactly one 1-bit here". Formally we e.g. have
 $(g,g,g):=\{(1,0,0),(0,1,0),(0,0,1)\}$ and 
 $$(g_1, {\bf 0},g_2,g_2,{\bf 1},g_1):=\big\{(1, {\bf 0},1,0,{\bf 1},0),(1, {\bf 0},0,1,{\bf 1},0),
 (0, {\bf 0},1,0,{\bf 1},1),(0, {\bf 0},0,1,{\bf 1},1) \big\}.$$

 {\bf 6.1} Let $G=(V,E)$ be any graph with a perfect matching $M\s E$. Then $\bigcup M=V$ and putting  $n:=|M|$ it follows that $|V|=2n$ and $\nu(G)=n$. If additionally $G=B$ is {\it bipartite} then each shore is a $n$-element vertex-cover\footnote{In particular $\tau(B)=n=\nu(B)$, so in {\it this} scenario $B$ being KE is obvious.}. 
  The sparser $B$ the more likely there will be mvc's beyond the two shores, and the more likely $MVC(B)$ will  be compressible with $g$-wildcards.
 As an extreme  case, if $B$ has $n$ connected components, all of them single edges, then $MVC(B)=(g_1,g_1,g_2,g_2,...,g_n,g_n)$ has cardinality $2^n$. 
 More interesting, take the connected bipartite graph $B_2=(V_2,E_2)$ in Figure 7a, which has the perfect 
matching $M_2=\big\{\{1,\ol{1}\},...,\{7,\ol{7}\}\big\}\s E_2$. 
It turns out (details later) that $MVC(B_2)=r_1\uplus\cdots \uplus r_8$ has $8+4+\cdots+1=26$ elements (Table 2). For instance, the fact that the {\it 01g-row} $r_1$ is contained in $ MVC(B_2)$ can be argued as follows. The vertices $\ol{1},\ol{3},\ol{4},\ol{7}$ cover all edges of $B_2$ except the three fat ones. In each of the fat edges one can pick exactly one point, indicated by the three $g$-wildcards of $r_1$. Uniting  $\{\ol{1},\ol{3},\ol{4},\ol{7}\}$ with $2^3$ choices of 3-sets (e.g. $\{2,\ol{5},6\}$ or $\{\ol{2},5,\ol{6}\}$) yields 8 distinct mvc's $Z$ of $B_2$, and they (as bitstrings) constitute the members of $r_1$.

\begin{tabular}{c|c|c|c|c|c|c|c||c|c|c|c|c|c|c||c} 
 & 1 & 2& 3 & 4 & 5 & 6 & 7 & $\ol{1}$ & $\ol{2}$ & $\ol{3}$ & $\ol{4}$ & $\ol{5}$& $\ol{6}$ & $\ol{7}$\\ \hline
 &   &  &   &   &   & &   &  &   &   &   &  & &\\ \hline
 $r_1=$ & 0 & $g_1$ & 0 & 0 & $g_2$ & $g_3$ & 0     & 1 & $g_1$ &1 & 1 & $g_2$ & $g_3$ & 1&  8\\ \hline
  $r_2=$ & 1 & $1$ & 0 & 0 & $g_2$ & $g_3$ & 0     & 0 & 0 &1 & 1 & $g_2$ & $g_3$ & 1&  4\\ \hline
  $r_3=$ & 1 & 1 & 1 & 0 & $g_1$ & $g_2$ & 0     & 0 & 0 &0 & 1 & $g_1$ & $g_2$ & 1&  4\\ \hline
   $r_4=$ & 1 & 1 & 1 & 1 & 1 & 1 & 0     & 0 & 0 &0 & 0 & 0 & 0 & 1&  1\\ \hline

   $r_5=$ & 0 & $g_1$ & 0 & 0 & $g_2$ & 1 & 1     & 1 & $g_1$ &1 & 1 & $g_2$ & 0 & 0&  4\\ \hline 
   $r_6=$ & 1 & 1 & 0 & 0 & $g$ & 1 & 1     & 0 & 0 &1 & 1 & $g$ & 0 & 0&  2\\ \hline
   $r_7=$ & 1 & 1 & 1 & 0 & $g$ & 1 & 1     & 0 & 0 &0 & 1 & $g$ & 0 & 0&  2\\ \hline
    $r_8=$ & 1 & 1 & 1 & 1 & 1 & 1 & 1     & 0 & 0 &0 & 0 & 0 & 0 & 0&  1\\ \hline  
 \end{tabular}

 {\sl Table 2: Packing the 26 minimum vertex-covers of $B_2$ into 01g-rows.}

 {\bf 6.2} The following concept will be crucial. A  pair $(Z,M)$ is a {\it rake} of a graph $G$ if $Z$ is a vertex-cover of $G$, and $M$ is a matching of $G$, and every $v\in Z$ belongs to some edge of $M$, and every edge of $M$ is incident with exactly one vertex in $Z$. In formulas: $|Z|=|M|$ and $|\{v,v'\}\cap Z|=1$ for all $\{v,v'\}\in M$.
 In particular, if a graph has a rake then it must be a KE graph. More can be said \cite[p.~223]{crama2011}:
 \begin{itemize}
   \item[(16)] {\it A graph is KE iff it has a rake.}
   \item[(17)] {\it A graph is KE iff for every maximum matching $M$ and every mvc $Z$ the pair $(Z,M)$ is a rake. }
 \end{itemize}
 
For instance, if $Z_1:=\{1,4,6\}$ and $M_1:=\{\{1,2\},\{3,4\},\{5,6\}\}$ then $(Z_1,M_1)$ is a rake of $G_1$ in Fig.2, and $B_2$ in Fig.7(a) has 26 rakes of type $(Z,M_2)$, where $M_2:=\{\{1,\ol{1}\},...,\{7,\ol{7}\}\}$.

{\bf 6.3} K\"onig-Egarv\'ary graphs with perfect matchings (such as $G_1$ or $B_2$) can be generated systematically as follows. With each pure 2-CNF $f=f(x_1,...,x_n)$ we associate\footnote{It corresponds to the concept of \cite[p.~210]{crama2011}, except that we are dealing with 2-CNFs as opposed to 2-DNFs. Furthermore, our 2-CNFs being pure, no loops occur in $G_f$.} the so called {\it matched graph} $G_f$ whose $2n$ vertices are all literals $x_1,\ol{x}_1,...,x_n,\ol{x}_n$. As to the edges of $G_f$, on the one hand they correspond to the clauses of $f$, on the other hand also\footnote{Some  clauses $x_i\vee\ol{x}_i$ may be present in $f$ as well but, being tautological, they can  be dropped. It merely matters  that all {\it  edges} $\{x_i,\ol{x}_i\}$ be present.}
$\{x_1,\ol{x}_1\},...,\{x_n,\ol{x}_n\}$ are postulated as edges.
Consequently they constitute a perfect matching  $M_f$ of $G_f$.
For instance, if the pure 2-CNF is
$$(18)\quad f_1:=(x_1\vee x_2)\wedge (\ol{x_1}\vee x_2)\wedge (\ol{x_1}\vee x_3)\wedge (x_2\vee x_3)$$
then its matched graph $G_{f_1}$ is rendered in Fig.8. The modelset $Mod(f_1)$ can be determined by inspection.  If $z\in Mod(f_1)$ has $z_1=0$, then $z_2=1$ (because of $x_1\vee x_2$). Hence $z_3$ can be both (because of $(0\vee 1)\wedge(\ol{0}\vee 1)\wedge(\ol{0}\vee z_3)\wedge(1\vee z_3)=1$). But if $z_1=1$, then again $z_2=1$ (because of $\ol{x}_1\vee x_2$), and so $z_3=1$ (because of $\ol{x}_1\vee x_3$). To summarize
$$Mod(f_1)=\{(0,1,0),(0,1,1),(1,1,1)\}.$$
Notice that the corresponding "literal models", i.e. 
$$(19)\quad \{\ol{x}_1,x_2,\ol{x}_3\},\ \{\ol{x}_1,x_2,x_3\},\ \{x_1,x_2,x_3\},$$ all happen to be vc's of $G_{f_1}$. Because $G_{f_1}$ has a matching of cardinality 3, these vc's are in fact mvc's. 

\begin{center}
\includegraphics[scale=1.1]{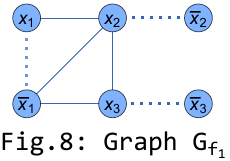}
\end{center}

{\bf 6.4} In order to generalize this in Theorem 3, we first elaborate the concept glimpsed in (19).
Given a Boolean formula $f:\{0,1\}^n\to\{0,1\}$ (that hence uses literals), a {\it literal\footnote{ While being a pun, let us embrace this name to capture an important concept that has been nameless so far.} model} of $f$ is a set of $n$ literals which all evaluate to 1. 
More precisely, each $(y_1,...,y_n)\in Mod(f)$ induces the literal model $\{x_i:\ y_i=1\}\cup \{\ol{x}_i:\ y_i=0\}$. Conversely, claiming that (say) $\{x_1,\ol{x}_2,\ol{x}_3,x_4\}$ is a literal model of $f:\{0,1\}^4\to\{0,1\}$ amounts to say that $(1,0,0,1)\in Mod(f)$.
We write $LitMod(f)$ for the set of all literal models of $f$. Hence $LitMod(f)\s {\cal P}(\{x_1,\ol{x}_1,...,x_n,\ol{x}_n\})$ is bijective to $Mod(f)\s \{0,1\}^n$ in straightforward ways. 
The next result matches Theorem 5.14 in \cite[p.~223]{crama2011}, which in turn is credited to \cite{Simeone}. We adapt Theorem 5.14 to our terminology.

{\bf Theorem 3:}
{\it  
\begin{itemize}
  \item[(a)] If $f$ is a pure and satisfiable 2-CNF, then $G_f$ is a KE graph with a perfect matching. Furthermore the literal models of $f$ coincide with the mvc's of $G_f$.
  \item[(b)] Each KE graph $H$ with a perfect matching is isomorphic to some matched graph $G_f$. 

\end{itemize}  }

{\it Proof.} As to (a), we previously stressed the crucial perfect matching of $G_f$.
But why is $G_f$ a KE graph? Because $f$  is satisfiable by assumption, models $z=(z_1,...,z_n)$ exist. If $Z$ is the corresponding literal model then by the satisfiability assumption $Z$ intersects all edges of $G_f$. So $Z$ is a vertex-cover. In fact
$(Z,M)$ is a rake of $G_f$ since $|Z\cap \{x_i,\ol{x}_i\}|=1$ for all $\{x_i,\ol{x}_i\}\in M$. Therefore $G_f$ is KE by (16). Since conversely each rake yields a unique literal model of $f$, there is a bijection between $LitMod(f)$ and $MVC(G_f)$.

As to (b), label the perfect matching $M$ of $H$ as $M=\big\{  \{x_1,\ol{x}_1\}, \ldots, \{x_n,\ol{x}_n\}\big\}$ and accordingly translate the edge-set of $H$ to a pure 2-CNF $f$. It follows that $H\simeq G_f$, {\it provided} $f$ is satisfiable. Why is $f$ satisfiable? Since $H$ is KE by assumption, by (17) any mvc $Z$ yields a rake $(Z,M)$.
Therefore for all $1\le i\le n$ exactly one literal in the edge  $\{x_i,\ol{x}_i\}$ belongs to $Z$, and at least one literal in every other edge of $H$ belongs to $Z$. Hence $Z$ is a literal model of $f$, and so $f$ is satisfiable. $\square$

Heeding the recipee in the proof of Thm.3(b) relabeling the perfect matching of $G_1$ in Figure 2 with literals happens to yield Figure 8. Hence $G_1\simeq G_{f_1}$. Looking at (19) we see that
$\ol{x}_2$ appears in no mvc.  Hence $\ol{x}_2\in core(G_{f_1})$. (Equivalently $3\in core(G_1)=A^*(G_1)$, see Fig.3.)

{\bf 6.5} Generally let $G=(V,E)$ be any KE graph with a $n$-element perfect matching $M$. Relabeling the vertices $v\in V$ admissibly with $2n$ literals (so each $\{v,v'\}\in M$  becomes some $\{x_i,\ol{x}_i\}$) triggers a pure 2-CNF $f=f(x_1,..,x_n)$ whose clauses correspond to the edges of $G$. If say $n=5$ and $y=(0,1,1,0,1)\in Mod(f)$ then, as seen, the literal model 
$\{\ol{x}_1,x_2,x_3,\ol{x}_4,x_5\}$ is a mvc of the (relabeled) graph $G$. It would be a bad idea to store millions of mvcs in this format, in particular when it means to forego a potential compression of $Mod(f)$! Instead of storing $\{\ol{x}_1,x_2,x_3,\ol{x}_4,x_5\}$ let us revamp it to
$y^*:=(0,1,1,0,1; 1,0,0,1,0)$, from which the former is readily recovered
(assuming that the components of all  $2n$-bitstrings are indexed by $x_1,...,x_5;\ol{x_1},...,\ol{x}_5)$. Moreover the compression of $Mod(f)$ can be salvaged as follows. Say $(2,0,2,1,2)\s Mod(f)$. Then the 8 revamped literal models corresponding to the 8 models in $(2,0,2,1,2)$ are collected in $(g_1,0,g_2,1,g_3;\ g_1,1,g_2,0,g_3)$, where we used the $g$-wildcards from 6.1. Generally, if $Mod(f)$ is compressed using $R$ many 012-rows then $MVC(G)$ can be compressed using $R$ double-length 01g-rows.

\section{Compressing $MVC(G)$ (and thus $MACL(G)$) for general KE graphs $G$ }

  We now tackle KE  graphs $G=(V,E)$ which {\it need not} have a perfect matching. Theorem 4 unravels the structure of $MVC(G)$ in this case. This leads to a output-polynomial (and potentially compressed) enumeration of $MVC(G)$, and thus of $MACL(G)$. The structure of $core(G)$ will be discussed as well.
  
{\bf 7.1}
In Theorem 4 below everything depends on a fixed chosen maximum matching, which (Section 2) can be found in time $O\big(\sqrt{|V|}\ |E|\big)$. 
 
\vspace{3mm}
{\bf Theorem 4: }{\it Let $G=(V,E)$ be a KE graph and $M$  a fixed maximum matching. 
\begin{itemize}
    \item[(a)]Putting $V_M:=\bigcup M$  the induced subgraph $H:=G[V_M]$ is
     KE with perfect matching $M$ and $MVC(G)\s MVC(H)$.
    
    \item[(b)] It holds that $A\s V_M$ and that the tail $T=T(G):=V\setminus V_M$ is an anticlique\\ contained in $ NO$.
   
    \item[(c)] If $NO_M$ is the set of siblings of $A$ with respect to $M$, then $NO=T\uplus NO_M$. 
 
    \item[(d)] It holds that $MVC(G)=\{A\uplus Z:\ Z\in MVC(G[F])\}$. Here the induced subgraph $G[F]$ is KE and has a perfect matching
    $M_F\s M$.
\end{itemize}
}

\begin{figure}[h]  
\centering
\includegraphics[scale=0.3]{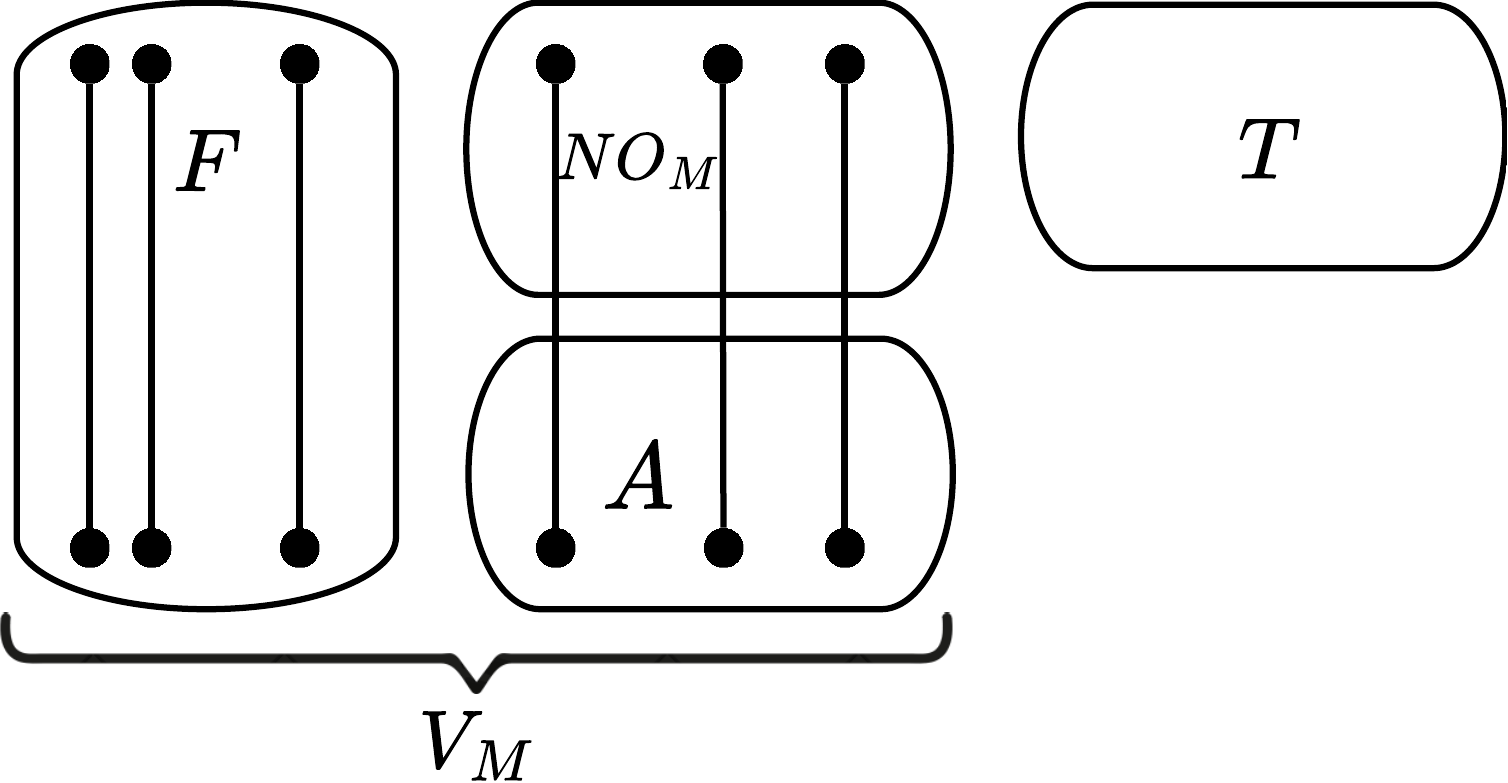}
\caption*{Figure 9: The A-F-NO Theorem for K\"onig-Egarv\'ary  graphs}
\end{figure}

As to Figure 9, recall from Section 3 the definition of $A,F,NO$. By duality statement (5) translates to:
$$(20)\ NH(NO)\s A,\ NH(F)\s A\uplus F,\ NH(A)\s A\uplus F\uplus NO$$
Furthermore  note that in Figure 9 only the edges of the maximum matching $M$ are indicated and that $NO=NO_M\uplus T$ by Thm.4(c). 

\vspace{3mm}
{\it Proof.}  
As to (a), since $G$ is KE and $M$ is a maximum matching, each mvc  $Z$ of $G$ yields a rake $(Z,M)$ of $G$ by (17). A fortiori $(Z,M)$ is a rake of $G[V_M]$. Hence $H:=G[V_M]$  is KE by (16). Clearly $M$ is a perfect matching of $H$ and
 $MVC(G)\s MVC(H)$. 

As to (b), no vertex in the {\it tail} $T:=T(G)$ belongs to any mvc $Z$ of $G$ (since $Z\s V_M$ by (a)). Hence $T\s NO$. From $A\cap T\s
A\cap NO=\es$ follows that $A\s V_M$. If there were adjacent vertices $v,v'\in T$ then $M\cup\{\{v,v'\}\}$ would be a matching larger than $M$. This contradiction shows that $T$ is an anticlique.

As to (c), each $v\in V_M$ has a unique {\it sibling} defined by $\{v,sib(v)\}\in M$. Since $(Z,M)$ is a rake for all $Z\in MVC(G)$, it holds for all $Z\in MVC(G)$ and all $v\in V_M$ that
$v\in Z\ \LRa\  sib(v)\not \in Z$.
In particular, it holds that 
$$(21)\quad (\forall v\in V_M)\quad v\in A\ \LRa\  sib(v) \in NO.$$ Hence if $NO_M:=\{sib(v'):\ v'\in A\}$
then $NO_M=NO\cap V_M$, and so $NO=T\uplus NO_M$.

As to (d), recall from Theorem 1' that $MVC(G)=\{Z\uplus A:\ Z\in MVC(G[F]) \}$ holds for {\it all} graphs $G$. As to the remaining claims,
clearly $M_F:=\big\{\{v,v'\}\in M:\{v,v'\}\s F\big\}$
is a perfect matching of $G[F]$ (see Fig.9). Take any $Z\in MVC(G)$. Each edge of $G$ contained in $F$ is touched not just by some (unpredictable) vertex $v\in Z$ but necessarily by some $v\in Z\cap F$. Hence $Z\cap F$ is a vertex-cover of $G[F]$, and from $|Z\cap F|=|M_F|$ follows that $G[F]$ is KE.
$\square$

 According to Theorem 4(d) for any KE graph $G$ the calculation of $MVC(G)$ can be reduced to the calculation of $MVC(G[F])$, where $G[F]$ is a certain KE graph  with a perfect matching. Hence $MVC(G[F])$ can be calculated via 2-CNFs by Theorem 3. Trouble is, in Theorem 4(d) one also needs to know $A(G)$. This (and more) is mitigated in the remainder of Section 7.

{\bf 7.2} Let $G=(V,E)$ be KE. If $M$ is a maximum matching and $V_M:=\bigcup M$ then the tail $T=V\setminus V_M$ is contained in $NO$ by Thm.4(b). The tail uniquely determines its {\it allmvc-shadow}, i.e. a certain set $A_S$ satisfying $A_S\s A:= A(G)$. For starters, $S_1:=NH(T)\s NH(NO)\s A$ by $(20)$. Hence $S_2:=sib(S_1)\s NO$ by $(21)$. Hence $S_3:=S_1\cup NH(S_2)\s A$ by $(20)$. Hence $S_4:=sib(S_3)\s NO$ by $(21)$. Hence
$S_5:=S_3\cup NH(S_4)\s A$ by $(20)$, and so forth.
The above process stabilizes as soon as $S_{i+2}=S_i=:A_S$ (where necessarily $i$ is odd). Because $NO_S:=sib(A_S)\s NO_M\s NO$ we call $NO_S$ the {\it nomvc-shadow} of $T$. Putting  $$(22)\quad V':=V_M\setminus (NO_S\uplus A_S)$$
(so $V'\supseteq F$) it holds by construction of $NO_S$ that no edges go from $NO_S$ to $V'$, and also none from $T$ to $V'$. Furthermore by construction of $A_S$ one has $A_S\cap\{u,v\}\neq\es$ for all edges $\{u,v\}\in E$ that satisfy $u,v\in T\uplus NO_S\uplus A_S$. Put another way:

\begin{itemize}
    \item[(23)] {\it $A_S$ is a vertex-cover of $G[ T\uplus NO_S\uplus A_S]$ and the only potential edges between $NO_S\uplus A_S$ and $V'$ are between  $A_S$ and $V'$. }
\end{itemize}

{\bf 7.3} For illustration consider the KE graph $G_4=(V,E)$ (Fig.10) which has a maximum matching $M_4$ consisting of the eight boldface vertical edges. Hence the three leftmost vertices constitute the tail $T_4$ of $G_4$. They are labeled with $"0"$ to stress the fact that they belong to {\it no} mvc of $G_4$. Likewise vertices that belong to {\it each} mvc will carry labels $"1"$ (with subscripts). Adopting the above notation the set $S_1=NH(T_4)$ contains the two vertices labeled $"1_1"$. Because set notation implies that $\{1_1,1_1\}=\{1_1\}$, we rather write $S_1=[1_1,1_1]$. Hence
$S_2=[0_2,0_2]$ and $S_3=[1_1,1_1,1_3,1_3]$
(one vertex $1_3$ is triggered by two $0_2$'s). Further $S_4=[0_2,0_2,0_4,0_4],\ S_5=[1_1,1_1,1_3,1_3,1_5],\ S_6=
[0_2,0_2,0_4,0_4,0_6],\ S_7=S_5=A_S$. The three dashed edges in Fig.10 are the ones that lead from $A_S$ to $V'$.
As stated in $(21)$, check that all neighbors of 0-labeled vertices are 1-labeled, and so $A_S$ is a vertex-cover of $G_4[T_4\uplus NO_S\uplus A_S]$. We mention in passing that the six rightmost vertices in Fig.10 induce a subgraph isomorphic to $G_1$.

\begin{center}
\includegraphics[scale=1.3]{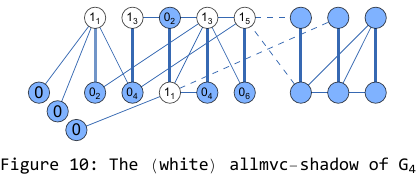}
\end{center}

{\bf 7.4} Theorem 5(i) below resembles Thm.4(d) in that both $G[F]$ and the larger graph $G[V']$ are KE graphs with a perfect matching. However, as opposed to $A(G)$ in Thm.4(d), the allmvc-shadow $A_S$ in Thm.5(i) can be calculated fast.

{\bf Theorem 5: }{\it Let $G=(V,E)\ (w:=|V|,m:=|E|)$ be a KE graph with allmvc-shadow  $A_S$ (induced by the tail $T(G)$). Then the following holds:
 \begin{itemize}
\item[(i)] $MVC[G]=\big\{A_S\uplus Z':\ Z'\in MVC(G[V'])\ \big\}$. Here $V'$ is defined in $(22)$ and is such that $G[V']$ is KE with a perfect matching.
\item[(ii)] The $N$ members of $MVC(G)$ can be calculated in time $O(Nw^3)$. Furthermore $MVC(G)$ can potentially be compressed by virtue of 01g-rows.

\item[(iii)]  $A(G)=A_S\uplus A(G[V'])$, and $A_S\neq\es$ when $G$ is connected.
\end{itemize}
}

{\it Proof.} As to (i), analogous to $G[F]$ in Thm.4(d) one shows that $G[V']$ is KE and has a perfect matching $M'\s M$ (so $M'$ contains those edges of $M$ which are disjoint from $A_S$). In order to
 show $\s$ in (i), take $Z\in MVC(G)$. Then $Z=A_S\uplus (Z\cap V')$ because of $Z\s V_M$ and $(22)$. 
 Furthermore $Z':=Z\cap V'$ belongs to $MVC(G[V'])$ because $|Z'|=|M'|$. As to $\supseteq$, let $Z'\in MVC(G[V'])$. Then $A_S\uplus Z'$ is a vertex-cover of $G$ by $(23)$, and a minimum one because $|A_S\uplus Z'|=|M|$. 
 
 As to (ii), find a maximum matching $M$ of $G$ in time $O(\sqrt{w} \ m)$ (Section 2) and calculate the two shadows $A_S,NO_S$ in even less time. Since $V'=V\setminus (A_S\uplus NO_S)$ and
 since $G[V']$ is KE and has a perfect matching $M'$, there is by Thm.3(b) a pure 2-CNF $f$ such that $G[V']\simeq G_f$ and  $MVC(G_f)=LitMod(f)$. By the definition of $N$ therefore $N=|MVC(G)|=|Mod(f)|$.
Since according to \cite[Thm.~1]{wild2025} one can enumerate $Mod(f)$ in time $O(Nw^3)$, it follows that $MVC(G)$ can be enumerated in time $O(\sqrt{w}\cdot m+Nw^3)=O(\sqrt{w}\cdot w^2+Nw^3)=O(Nw^3)$ (viewing that $N\neq 0$). As to the second claim in (ii), in 6.5 we argued that the compression of $Mod(f)$ via 012-rows triggers a compression of
$MVC(G[V'])$ with equally many 01g-rows. Upon extending each 01g-row  with $|A_S|$ many 1-bits one obtains, by virtue of (i), a compression of $MVC(G)$.

 As to (iii), the first part follows at once from (i). Looking at 7.2 the connectedness of $G$ implies that $S_1=NH(T)\neq\es$. By definition of $A_S$ this implies that $A_S\neq\es$.
 $\square$

{\bf 7.5} Let us translate Theorem 5 from a mvc-framework to a macl-framework. We begin by recalling old terminology and introducing new one (in all four cases "shadow" means "shadow of the tail $T^*(G):=T(G)$") :
\begin{itemize}
    \item {\it allmacl-shadow} $A_S^*:=NO_S$ (= our previous nomvc shadow)
    \item {\it nomacl-shadow} $NO_S^*:=A_S$ (= our previous allmvc shadow)
\end{itemize}

{\bf Theorem 6: } {\it For each KE graph $G=(V,E)$ with fixed maximum matching $M\s E$ (and the ensuing parameters defined above) the following holds.
\begin{itemize}
    \item [(i)] $MACL(G)=\big\{T^*(G)\uplus A_S^*(G)\uplus Y':\ Y'\in MACL(G[V'])\big\},$
where $G[V']$ is an induced KE graph with a perfect matching.
\item[(ii)]  The $N$ members of $MACL(G)$ can be calculated in time $O(Nw^3)$.
\item[(iii)] $A^*(G)=T^*(G)\uplus A_S^*(G)\uplus A^*(G[V'])$, and so $A^*(G)=core(G)$ can\\ be calculated in time $O(w^3)$.
\end{itemize} }

{\it Proof.} As to (i), it follows from $(22)$ and the above definitions that
$$(24)\quad V_M\ =\ NO_S\uplus A_S\uplus V'\ =\ A^*_S\uplus NO^*_S\uplus V'$$
Recalling that in each graph the macls are the complements of the mvcs it holds that
$$MACL(G)=\big\{V\setminus Z:\ Z\in MVC(G)\big\}\stackrel{T.4(b)}{=\!=}\big\{(T(G)\uplus V_M)\setminus Z:\ Z\in MVC(G)\big\}$$
$$\stackrel{T.4(a)}{=\!=}\big\{T(G)\uplus (V_M\setminus Z):\ Z\in MVC(G)\big\}
\stackrel{T.5(i)}{=\!=}\big\{T(G)\uplus (V_M\setminus (A_S\uplus Z')):\ Z'\in MVC(G[V'])\big\}$$
It follows from from (24)  that the sets $V_M\setminus(A_S\uplus Z')$ are exactly\footnote{A visualization is obtained by adapting Figure 4 from universe $V$ to universe $V_M$: Relabel $F$ by $V'$; $Z$ by $Z'$; $Y$ by $Y'$; $NO$ by $NO_S$; $A$ by $A_S$. } the sets $A_S^*\uplus Y'$ when $Y'$ ranges over  $MACL(G[V'])$. 

As to (ii),  $MACL(G[V'])=\{V'\setminus Z':\ Z'\in MVC(G[V'])$ (Section 2). Recall from the proof of Thm.5(ii) that $MVC(G[V'])$ can be compressed with at most $N:=|MVC(G[V'])|$ many 01g-rows $r_i$ in time $O(Nw^3)$. Putting 
$r_i^c:=\{V'\setminus Z':\ Z'\in r_i\}$ it follows that $MACL(G[V'])=r_1^c\uplus\cdots\uplus r_N^c$. Each set-system $r_i^c$ can again be written as 01g-row, namely as the 01g-row obtained from $r_i$ by switching 0's and 1's, while leaving all g-symbols unchanged. It follows from (i) that $MACL(G)$ can be compressed with at most $N$ many 01g-rows in time $O(Nw^3)$.
(If this is too confusing, by Thm.5(ii) enumerate all mvcs $Z$ one-by-one in output polynomial time. Then the macls are the sets $V\setminus Z$.)

The identity in (iii) follows at once from the identity in (i). Hence the calculation of $core(G)$ boils down to the calculation of $core(G[V'])$. Let us generally handle $core(H)$ for any KE graph $H$ with perfect $n$-element matching. Referring to Thm.3 (a),(b) let $f$ be the 2-CNF for which $LitMod(f)$ is bijective to $MVC(H)$.
The elements of $core(H)$ are the elements of $H$ which appear in no mvc of $H$. These elements correspond to the literals which appear in no literal model. 
The literal $x_k$ appears in no literal model iff the (non-pure) 2-CNF $x_k\wedge f$ is insatisfiable. Likewise the literal $\ol{x}_k$ appears in no literal model iff the (non-pure) 2-CNF $\ol{x}_k\wedge f$ is insatisfiable. In view of Sec. 4 these $2n\le w=|V|$ satisfiability tests can be carried out in time $O(2n(n+m))=O(w(w+w^2))=O(w^3)$.
$\square$

\section{Specialization to bipartite graphs}

We first show (8.1 to 8.3) how Theorem 6 can be trimmed in the bipartite case, both conceptually and computation-wise. Then we point out connections to the so called Dulmage-Mendelson decomposition of a bipartite graph.

{\bf 8.1} In Section 3 we stated (I) that the $A^*-F^*-NO^*$ Theorem (=Theorem 1) is the more relevant the easier $A^*$ can be determined and the more inviting the structure of $G[F^*]$ is. Also recall from (8) that (II) for each graph $G$ each $Y\in MACL(G[F^*])$ can be matched to a disjoint $Y_0\in MACL(G[F^*])$. Finally (III), whenever $G$ is connected and satisfies $\alpha(G)>\nu(G)$ then
$|core(G)|\ge 1+\alpha(G)-\nu(G)$.

How can (I),(II),(III) be refined for KE graphs $G$? As to (I), using 2-CNFs one can determine $A^*(G)$ in time $O(w^3)$ by Thm.6(iii), and $MACL(G)$ in output-polynomial time $O(Nw^3)$ by Thm.6(ii). As to (II), this is property follows at once from the fact that the {\it whole} of $G[F^*]$ has a perfect matching (by Thm.4(d) and $F^*=F$). As to (III), when $G$ is connected then $A^*_S\neq\es$ by the dual of Thm.5(iii). In this case Thm.6(iii) implies
$$|core(G)|=|T^*(G)|+|A_S^*(G)|+|core(G[V'])|\ge (\alpha(G)-\nu(G))+1+0.$$

For bipartite graphs $G$ it gets better still. Namely, by Theorem 7 below,  the core $A^*(G)$ in Thm.6(iii) boils down to the readily computed set $A^*_S(G)\uplus T^*(G)$, and $MACL(G[F^*])$ can be handled by 2-CNFs with beneficial {\it mixed} clauses (details in 8.2).

{\bf Theorem 7: }{\it For each bipartite graph $B=(V,E)$ it holds that\\
$A^*(B)=T^*(B)\uplus A_S^*(B)$. Therefore $A^*(B)=core(B)$ can be calculated in time $O(w^{2.5})$. }

{\it Proof.} Recall that $A^*_S(B)$ is the allmacl-shadow of the tail $T^*(B)$. By Theorem 6 (iii) it suffices to show that $core(B[V'])=\es$. Indeed, since the induced graph $H:=B[V']$ stays bipartite and has a perfect matching, the two shores of $H$ are {\it disjoint} macls, which implies that $core(H)=\es$.
As opposed to Thm.6(iii) testing the satisfiability of 2-CNFs is no longer required. However, due to finding a maximum matching the overall cost only slightly reduces from $O(w^3)$ to $O(w^{2.5})$. $\square$

\begin{figure}[h]   
\centering
\includegraphics[scale=0.33]{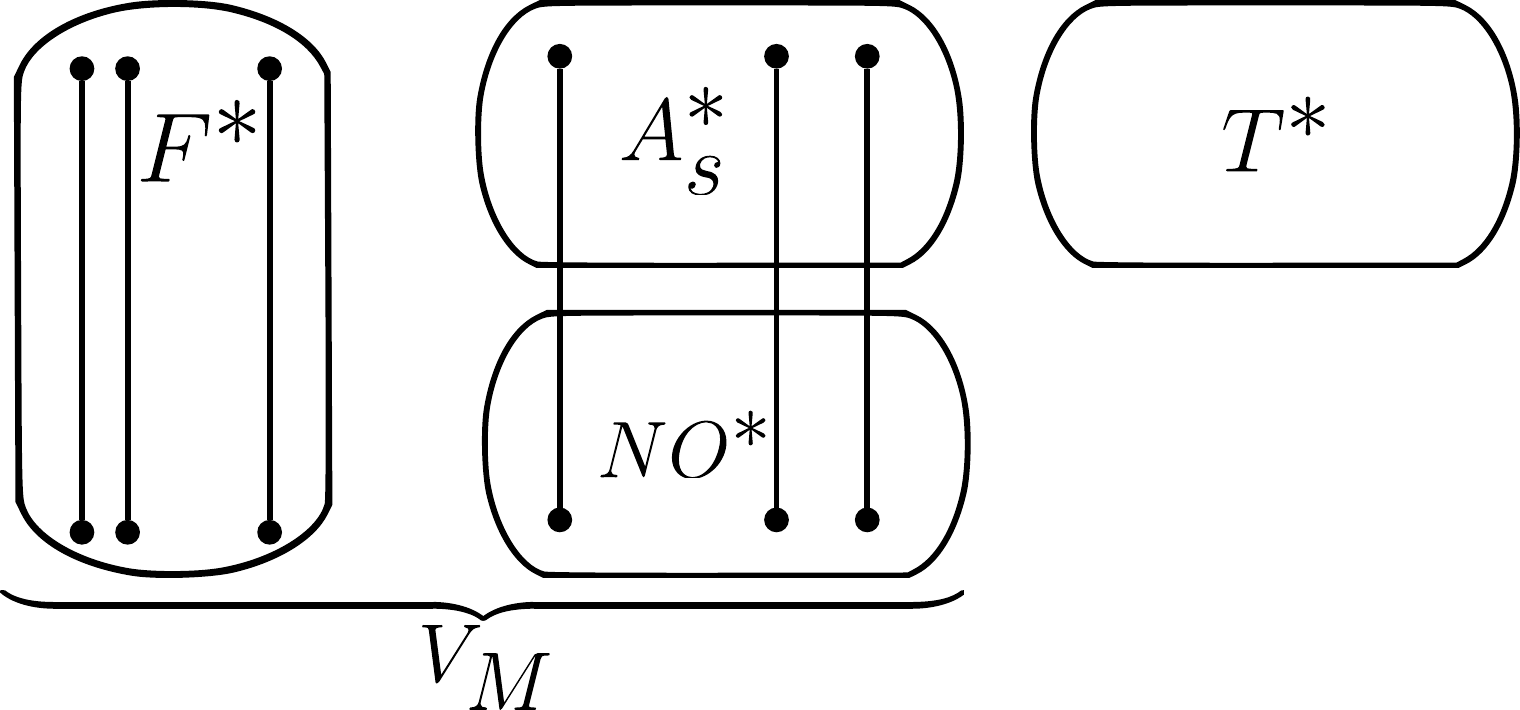}
\caption*{Figure 11: The A*-F*-NO* Theorem for bipartite graphs}
\end{figure}

By duality the visualization of the  $A^*-F^*-NO^*$ Theorem for KE graphs is obtained from Figure 9 upon changing the acronyms as follows: $F,T $ become $F^*,T^*$; $NO_M$ becomes $A_M^*$; $A$ becomes $NO^*$. In the bipartite case we have $A_M^*=A_S^*$ (Theorem 7), and hence obtain Figure 11.

{\bf 8.2} Let us unveil how in 5.5 the principal structure coupled to $B_2$ was calculated.
To begin with, relabeling its perfect matching $M_2$ with literals leads to Figure 12a.
As for all bipartite graphs with perfect matching (such as all $B[F^*]$), the triggered pure 2-CNF $f_2$  has exclusively
 {\it mixed} clauses (as opposed to $f_1$ in (16)). So e.g. edge $\{\ol{x}_4,x_6\}$ yields the clause $\ol{x}_4\vee x_6$.
 Rather than immediately feeding $f_2$ to the algorithm in \cite{wild2025} (which was used for KE graphs in the proof of Thm 5(ii)) let us process $f_2$ further. 
 Namely, e.g. $\ol{x}_4\vee x_6$ can be written as  $x_4\ra x_6$. 
Directing all edges of $G_{f_2}$ in this way (and dropping the edges in $M_2$) yields the digraph $DG(f_2)$ in Fig.12b
Since $DG(f_2)$ happens to have no directed cycles, it yields a poset $P(f_2)$ where the order relation $\ge$ is defined as
$$x_i\ge x_j\ \LRa\ i=j\ or\ there\ is\ a\ directed\ path\ x_i\ra x_k\ra\cdots\ra x_j.$$
Figure 12c is the ensuing poset diagram. Hence $u=(u_1,...,u_7)\in\{0,1\}^7$ is a model of $f_2$  iff $u_i\ra u_j$ evaluates to 1 for all arcs $x_i\ra x_j$ in Fig.12b. Because these arcs $u_i\ra u_j$ are never of type $1\ra 0$, a moment's thought shows that $u\in\{0,1\}^7$ is a model of $f_2$ iff
the set of literals $\{x_i: \ u_i=1\}$ is an ideal of $P(f_2)$. Notice that $P(f_2)$ is isomorphic to the poset in Fig.7(b).

\begin{center}
\includegraphics[scale=1.1]{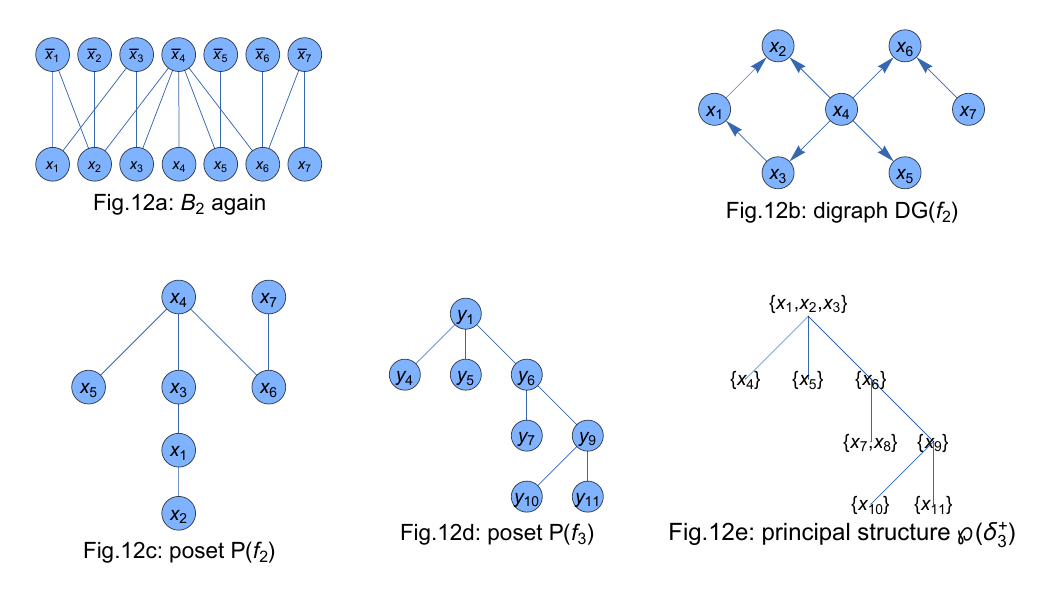}
\end{center}

In each poset the set of its  minimal elements constitutes an  ideal $K^{min}$ with the perk that each subset of $K^{min}$ is an ideal as well. For the poset $P(f_2)$ in Fig.12c one has $K^{min}:=\{x_2,x_5,x_6\}$.  Written as literal model this  becomes the mvc
 $\{\ol{x}_1,x_2,\ol{x}_3,\ol{x}_4,x_5,x_6,\ol{x}_7\}$ of $B_2$, which however (see 6.5) we rather store as
$Z:=(0,1,0,0,1,1,0;\ 1,0,1,1,0,0,1)$. It is one of 8 mvcs
 comprised in the 01g-row $r_1$ of Table 2. What about the other 7?
We have $(0,2,0,0,2,2,0)\s Mod(f_2)$ because all 8 subsets of $K_{min}=\{x_2,x_5,x_6\}$ remain ideals of $(P(f_2),\le)$. According to 6.5 the 012-row $(0,2,0,0,2,2,0)$ gets stretched
to the 01g-row $(0,g_1,0,0,g_2,g_3,0;\ 1,g_1,1,1,g_2,g_3,1)$, which equals $r_1$ in Table 2.
 Vice versa, let us e.g. look at the 01g-row $r_6$ in Table 2 and ask what ideals of $P(f_2)$ are encoded by the two members of $r_6$. One ideal is 
$\{x_1,x_2,x_6,x_7\}$ (check with Fig.12C that this is indeed an ideal). The other ideal is $\{x_1,x_2,x_6,x_7,x_5\}$.

{\bf 8.3} Here comes a larger bipartite graph $B_3$ whose diagram we omit and whose perfect matching $M_3$ has been relabeled with literals already, say $M_3=\big\{\{x_1,\ol{x}_1\},...,\{x_{11},\ol{x}_{11}\}\big\}$. As 8.1 the triggered pure 2-CNF $f_3$ could be
$$ f_3:=(x_1\ra x_2)\wedge(x_2\ra x_3)\wedge(x_3\ra x_1)\wedge(x_1\ra x_4)\wedge(x_3\ra x_5)\wedge(x_3\ra x_6)$$
$$\wedge\ (x_6\ra x_8)\wedge(x_6\ra x_9)\wedge(x_7\ra x_8)\wedge(x_8\ra x_7)\wedge(x_9\ra x_{10})\wedge(x_9\ra x_{11})$$
While again $MVC(B_3)=LitMod(f_3)$, for the time being let us focus on $Mod(f_3)$.
As in 8.1 the constituents $x_i\ra x_j$ are the arcs of the digraph $DG(f_3)$.
But different from 8.1 we now confront two directed cycles in $DG(f_3)$, i.e. $x_1\ra x_2\ra x_3\ra x_1$ and $x_7\ra x_8\ra x_7$.
It helps that they  happen to coincide with the strong components in which they are contained\footnote{Recall that a strong component can contain {\it several} directed cycles.}. Having condensed the strong components $\{x_1,x_2, x_3\}$ and $\{x_7,x_8\}$ to $y_1$ and $y_7$ respectively, and upon identifying all other $x_i$ with $y_i$, the 2-CNF $f_3$ becomes a new (but much related) 2-CNF:
$$ f_3':=(y_1\ra y_4)\wedge(y_1\ra y_5)\wedge(y_1\ra y_6)\wedge(y_6\ra y_7)\wedge(y_6\ra y_9)\wedge(y_9\ra y_{10})\wedge(y_9\ra y_{11}).$$
Specifically, it holds that 
$$Mod(f_3)=\big\{(y_1,y_1,y_1,\ y_4,y_5,y_6,\ y_7,y_7,\ y_9,y_{10},y_{11}):\ (y_1,y_4,y_5,y_6,y_7,y_9,y_{10},y_{11})\in Mod(f_3')\big\}$$
See Fig.12d. If we imagine that $B_3$ was handled as in Section 5 via some clumsy network with  submodular function $\delta_3^+$ then Fig.12e would be the ensuing principal structure.

It just so happens that $f_3'$ by distributivity coincides with (10) in 4.2. Therefore Table 1 compresses $Mod(f_3')$ by way of 012ab-rows. This yields $|Mod(f_3')|$ at once and therefore also $$(25)\quad |MVC(B_3)|=|LitMod(f_3)|=|Mod(f_3)|=|Mod(f_3')|.$$

Wrapping up, for general KE graphs $G$ Thm.5(ii) provided a compression of $MVC(G)$ by virtue of 01g-rows. The latter were obtained by "stretching" the 012-rows that constitute
the compression of $Mod(f)$ (where $f$ is a certain pure 2-CNF). As seen in 8.2 and 8.3, when
$G=B$ is bipartite then $Mod(f)$ can even be compressed by 012ab-rows.
For some purposes, e.g. computing the mere cardinality as in (25) or random sampling\footnote{To fix ideas, say $r=(a_1,b_1,b_1,0,2,1,a_2,b_2)$ is a random 012ab-row constituting $Mod(f)$. Pick a random bitstring in $r$, say $(0,1,0,0,0,1,1,1)$. It triggers the mvc $(0,1,0,0,0,1,1,1; 1,0,1,1,1,0,0,0)$.}, one can stick with $Mod(f)$. If one insists on a compression of the (equicardinal) family $MVC(B)$, it will necessarily be less compact than the compression of $Mod(f)$. Namely, each 012ab-row $r$ expands in obvious ways as disjoint union of $2^t$ many 012-rows $r_i$ (where $t$ is the number of (a,b)-wildcards in $r$), and in turn each $r_i$ contributes a 01g-row to $MVC(B)$.

{\bf 8.4} Here comes the Dulmage-Mendelson decomposition. To begin with, in a bipartite graph $B$ two vertices of $NO^*(B)$ may well be adjacent\footnote{Thus in Figure 11 the sets of upper vertices of $M$, respectively lower vertices of $M$, must not be mistaken as being part of the upper respectively lower shore.} but then evidently belong to different shores. While two vertices in $A^*(B)$ cannot be adjacent, it will nevertheless be crucial in 8.4 to split $A^*(B)$  shore-wise.

\begin{center}
\includegraphics[scale=1.1]{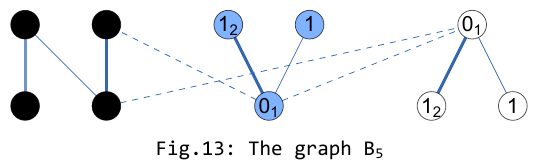}
\end{center}

To fix ideas consider the bipartite graph $B_5=(V_5,E_5)$ with boldface matching $M_5$ in Figure 13. Dually to Fig.10 we label the two vertices of the tail $T^*$ with "1" to stress that they belong to {\it all} macls, and then proceed to calculate (akin to 7.3) the allmacl-shadow $A_S^*$ and the equicardinal nomacl-shadow $NO_S^*$. One checks that $NO^*_S\ "="\ \{0_1,0_1\}=NO^*,\ A^*_S=\{1_2,1_2\}$, and so $A^*=T^*\uplus A^*_S=\{1,1,1_2,1_2\}$.

For general bipartite graphs $B$ with fixed maximum matching $M$ let us rephrase the construction of $NO^*(B)$ and $A^*(B)$ in terms of {\it $M$-alternating paths}, a concept that arises in many scenarios \cite{schrijver2003},\cite{LP},\cite{HV}. Namely, $NO^*(B)$ is the set of vertices that can be reached by an odd length $M$-alternating path that starts in $T^*(B)$. Likewise $A^*(B)$ is the set of vertices that can be reached by an even length $M$-alternating path that starts in $T^*(B)$. (In particular, for the vertices in $T^*(B)\s A^*(B)$ this path has length 0.) Let us put
$$Even=Even(B):=A^*(B),\quad Odd=Odd(B):=NO^*(B),\quad Unr=Unr(B):=F^*(B).$$
The acronym $Unr(B)$ reflects the fact that the vertices in this set are {\it unreachable} by all $M$-alternating paths. If $B=(V,E)$ is connected then $V=S_{up}\uplus S_{\ell o}$ where 
$S_{up}$ and $S_{\ell o}$ are the unique shores of "upper" and "lower" vertices. We put
$$Even_{up}:=Even\cap S_{up},\ Even_{\ell o}:=Even\cap S_{\ell o},
\ Odd_{up}:=Odd\cap S_{up},\ Odd_{\ell o}:=Odd\cap S_{\ell o}$$
and define these induced subgraphs of $B$:
$$B_{wc}:=B[Unr^*],\quad B_{oc}:=B[Even_{up}\uplus Odd_{\ell o}],\quad B_{uc}:=B[Odd_{up}\uplus Even_{\ell o}]$$

For $B=B_5$ these three graphs are colored black, blue and white in Figure 13. This {\it Dulmage-Mendelson (DM) decomposition} of any  bipartite $B$ is independent of the particular choice of $M$ and has many properties. To begin with, if $B$ is connected, then so are the components
$B_{wc},B_{oc},B_{uc}$. From $Even(B)=A^*(B)$ and (5) follows that there are no edges from the top of $B_{oc}$ to $B_{wc},B_{uc}$; nor from the bottom of 
$B_{uc}$ to $B_{wc},B_{oc}$. Apart from these exceptions, there may be edges between the components (and in Fig.13 they are dashed). But such edges never belong to a maximum matching.
In fact, each maximum matching $M_0$ of $B$ arises as follows. It consists of a perfect\footnote{That $B_{wc}$ {\it has} perfect matchings can be argued from
$B_{wc}=B[F]$ and Thm.4(d).} matching $M_{wc}$ of $B_{wc}$, a  matching $M_{oc}$ from $Odd_{\ell o}$ to $Even_{up}$ (all vertices of $Odd_{\ell o}$ are covered, at least one vertex of $Even_{up}$ is not), and
a  matching $M_{uc}$ from $Odd_{up}$ to $Even_{\ell o}$ (all vertices of $Odd_{up}$ are covered, at least one vertex of $Even_{\ell o}$ is not). 

{\bf 8.5} Some more facts about the DM decomposition of $B=(V,E)$. As seen, it arises from the $A^*-F^*-NO^*$ decomposition of $V$ by splitting $A^*$ and $NO^*$ according to shores and rearranging the parts.
It sheds light on the set of all maximum matchings but does not tell us how to enumerate them all (but for forests an elegant answer is given in Sec.9). In contrast, the $A^*-F^*-NO^*$ decomposition sheds light on $MACL(B)$, and we exploited that to also enumerate the whole of $MACL(B)$. 

The main application of the DM-decomposition is to systems of equations. Namely, the coupled bipartite graph $B$ has one vertex per equation $y_i$ and one vertex per unknown $x_j$, and there is an edge between $y_i$ and $x_j$ if $x_j$ appears in $y_i$. For instance, let $\{y_1,...,y_5\}$ and $\{x_1,...,x_5\}$ be the upper and lower shores in Figure 13. Then restricted to $B_{uc}=\{y_1,x_4,x_5\}$ the system is {\it underconstrained} (more variables than equations), restricted to $B_{oc}=\{y_3,y_4,x_3\}$ the system is {\it overconstrained} (less variables than equations), and restricted to $B_{wc}=\{y_1,y_2,x_1,x_2\}$ the system is {\it wellconstrained} (more variables than equations). We recommend \cite{AJM} as a gentle introduction to these matters. The DM-decomposition can also be viewed as principal
partition with respect to some submodular function \cite{Fu2008}, and as a special case of the Edmonds-Gallai-decomposition that applies to arbitrary graphs \cite{LP}. 

{\bf 8.6} The enumeration of all macls of a  bipartite (not  a general KE) graph $B$ was previously claimed\footnote{Some results are not fully proved, e.g. their crucial reference [6] is not accessible by Google Scholar.} in \cite[Thm.~2]{kashiwabara1992}. Let us point out differences and similarities.
Like us \cite[Thm.~2]{kashiwabara1992} hinges on a initial fixed maximum matching $M$. In a "First reduction" $B$ gives way to some subgraph $B'$ with  perfect matching that corresponds to our subgraph $B[V']$. (While $M$-alternating paths are used, there is no mention of Dulmage-Mendelson.) In a "Second Reduction" the macls of $B'$ are enumerated by considering, like us, an auxiliary digraph and its strongly connected components. Unfortunately the authors  \cite{kashiwabara1992} never proceed to the more comfortable level of posets and ideals, let alone to the  compressed enumeration of all ideals.

\section{Specialization to forests}

We give a brief proof for the fact that all trees have the KE property. This leads to a natural way to enumerate all maximum matchings of a tree. We close with comments on the article [Z] of Jennifer Zito.

{\bf 9.1} By definition forests $Tr$  are  cycle-free bipartite graphs $B$, and trees are connected forests. (The acronym Tr does not insinuate that we are dealing with trees. Even if one prefers trees, inductive arguments may require to consider forests.) Proving the KE property for forests $Tr=(V,E)$ is much easier\footnote{Surprisingly, I did not find the upcoming argument in the literature.} than for bipartite graphs.
Namely, we induct on the cardinality of $|V|$ and can assume that $E\neq\es$. Pick any leaf $v\in V$ and let $u$ be the unique vertex adjacent to $v$; such an edge $\{u,v\}$ will be called a {\it leaf-edge}.
By induction the forest $Tr[V\setminus \{u,v\}]$ has a matching $M'$ and an equicardinal vertex-cover $Z'$. It
follows that $M := M'\cup \{\{u, v\}\}$ is a matching of $Tr$ and that $Z'\cup\{u\}$ is an equicardinal vertex-cover of $Tr$. Hence $Tr$ is KE.

{\bf 9.2} As a consequence one obtains a fast algorithm to find a maximum matching in a forest $Tr$: Keep on pulling leaf-edges one after the other until one arrives at a "final forest" in the sense that all its connected components are isolated vertices or edges. The latter edges constitute a maximum matching of $Tr$. To illustrate, consider the tree $Tr_6$ in Fig.14a. Pulling from it (say) the leaf-edge $1$ and then $9$ yields the forest in Fig.14b. Pulling from it the leaf edge $4$ (pulling the leaf-edge $6$ has the same effect) yields the final forest in Fig.14c. Hence $\{1,4,6,9\}$ is a maximum matching of $Tr_6$. The presented method extends in straightforward ways (details omitted)  to a depth first search algorithm that enumerates in output-linear time all maximum matchings of a forest.

\begin{figure}[h]   
\centering
\includegraphics[scale=0.96]{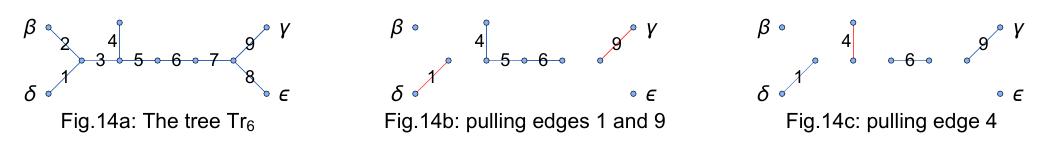}
\caption*{Figure 14: Getting a maximum matching of a tree}
\end{figure}

As to {\it perfect} matchings, even more can be said (this is e.g discussed in the Math stack exchange):
$$(26)\ Each\ forest\ has\ at\ most\ one\ perfect\ matching. $$
Here comes yet another proof. Let $B$ be bipartite with two distinct perfect matchings
$M$ and $\ol{M}$. It suffices to show that there is a cycle in $B$. Thus let $P$ be an alternating path of maximum length $t$, in the sense that its edges alternate between belonging to $M$ and $\ol{M}$. From $M\neq \ol{M}$ follows $t>1$. Let $e_1,e_2,...,e_t$ be the edges of $P$, from beginning to end. Without loss of generality $e_t=\{v,v'\}\in M$ and $v'$ is the endpoint of $P$ (so $v\in e_{t-1}\cap e_t$).
Since $\ol{M}$ is a perfect matching, there is exactly one edge $e$ of type $\{v',v''\}$ in $\ol{M}$. Since $P$ has maximum length, $v''$ is one of the vertices of $P$. If $v''\neq v$ this triggers a cycle in $B$. But why is $v''\neq v$? Because $v''=v$ implies $e_t=e\in \ol{M}$, which in view of $e_{t-1}\in\ol{M}$ is impossible. $\square$

Recall that for each KE graph $G$ the subgraph $G[F^*]$ induced by the flexible vertices has some perfect matching $M_{F^*}$. It follows from (26) that $M_{F^*}$ is uniquely determined when $G$ is a forest.


{\bf 9.3} For a tree $Tr=(V,E)$  and  $e\in E$ put $Tr^{-}:=(V,E\setminus\{e\})$. In \cite{zito1991} the edge $e$ is called $\alpha$-{\it critical} if $\alpha(Tr^{-})>\alpha(Tr)$. For instance $\alpha(Tr_6)=10-\tau(Tr_6)=10-\nu(Tr_6)=6$, from which it follows that edge $6$ is $\alpha$-critical, but edge $1$ is not.
Section 3 and   Theorem 4 of our article were inspired by \cite{zito1991} which for  trees $Tr=(V,E)$ introduced the partition $V=A^*\uplus F^*\uplus NO^*$. 

In 9.2 we argued that $Tr[F^*]$ has a unique  maximum matching $M_{F^*}$. In fact
$M_{F^*}$ coincides (Thm.2 in \cite{zito1991}) with the set of all $\alpha$-critical edges.


\section{Declaration about the use of AI}

Upon submitting the 45-page article \cite{wild2026} to the arXiv the author asked the AI Claude 5 to review it. The feedback was amazing. Leaving the positive comments aside, let us dwell upon the constructive criticism. To begin with, I heeded the advice to split the article. One half will be attended to in the future, the other half developped to the article in front of you. Here come the three main aspects that were sparked by Claude 5 (the second required the most extra input of mine). First, the incorporation of material from \cite{BML} into Section 3. Second,  the link to directed networks \cite[Lemma 4]{PQ} in Section 5. Third, the Dulmage-Mendelson decomposition of bipartite graphs in Section 8. 

\section{Further declarations}

 This research received no external funding. There is no conflict of interest. Data sharing is not applicable to this article as no datasets were generated or analyzed during the current study.

\section{References}

\printbibliography[heading=none]

@book{crama2011,
  author    = {Crama, Yves and Hammer, Peter L.},
  title     = {Boolean Functions: Theory, Algorithms, and Applications},
  series    = {Encyclopedia of Mathematics and its Applications},
  number    = {142},
  publisher = {Cambridge University Press},
  year      = {2011},
}

@article{Xiao,
  author  = {M. Xiao and H. Nagamochi},
  title   = {Exact algorithms for maximum independent set},
  journal = {Information and Computation},
  volume  = {255},
  year    = {2017},
  pages   = {126--146},
}

@article{Micali,
  author  = {S. Micali and  V. Vazirani},
  title   = { An $O(\sqrt{|V|} |E|)$ algorithm for finding a maximum matching in general graphs},
  journal = {21st Annual Symposium on Foundations of Computer Science, Syracuse, NY, USA.},
  volume  = {},
  number  = {},
  year    = {1980},
  pages   = {17--27},
}

@article{AJM,
  author  = {S. Ait-Aoudia and R. Jegou and D. Michelucci},
  title   = {Reduction of constraint systems},
  journal = {arXiv:1405.6131 },
  volume  = {},
  year    = {},
  pages   = {},
}

@article{Simeone,
  author  = {B. Simeone},
  title   = {Consistency of quadratic Boolean equations and the K\"onig-Egarv\'ary property for graps},
  journal = {Annals of Discrete Mathematics},
  volume  = {25},
  year    = {1985},
  pages   = {281--290},
}

@article{gainerdewar2017,
  author  = {Gainer-Dewar, A. and Vera-Licona, P.},
  title   = {The minimal hitting set generation problem: algorithms and computation},
  journal = {SIAM Journal on Discrete Mathematics},
  volume  = {31},
  year    = {2017},
  pages   = {63--100},
}

@article{Wild2,
  author  = {M. Wild},
  title   = {Compression with wildcards: All exact or all minimal hitting sets},
  journal = {Open Mathematics},
   volume  = {21},
  year    = {2023},
 
}

@article{kashiwabara1992,
  author  = {Kashiwabara, T. and Masuda, S. and Nakajima, K. and Fujisawa, T.},
  title   = {Generation of maximum independent sets of a bipartite graph and maximum cliques of a circular-arc graph},
  journal = {Journal of Algorithms},
  volume  = {13},
  year    = {1992},
  pages   = {161--174},
}

@book{schrijver2003,
  author    = {Schrijver, A.},
  title     = {Combinatorial Optimization},
  publisher = {Springer},
  year      = {2003},
}

@article{wild2025,
  author  = {Wild, M.},
  title   = {Compression with wildcards: All models of a {Boolean} 2-{CNF}},
  journal = {Contemporary Mathematics},
  volume  = {6},
  number  = {6},
  year    = {2025},
  pages   = {7795--7816},
}

@article{wild2026,
  author     = {M. Wild},
  title      = {Compression with wildcards: All, or all maximum, anticliques of a graph},
  journal = {arXiv},
  year = {2026},
}

@article{zito1991,
  author  = {Zito, J.},
  title   = {The structure and maximum number of maximum independent sets in trees},
  journal = {Journal of Graph Theory},
  volume  = {15},
  year    = {1991},
  pages   = {207--221},
}

@article{PY,
  author  = {S Pettie and L Yin},
  title   = { The Structure of Minimum Vertex Cuts},
  journal = {arXiv:2102.06805v1},
  volume  = {},
  year    = {2021},
  pages   = {},
}

@article{LP,
  author  = {L. Lovasz and D. Plummer},
  title   = {Matching Theory},
  journal = {Chelsea AMS Publishing},
  volume  = {},
  year    = {1986},
  pages   = {},
}

@article{Fu2008,
  author  = {S. Fujishige},
  title   = {Theory of Principal Partitions Revisited},
  journal = {Trends in Combinatorial Optimization: Bonn 2008, 2009 - Springer},
  volume  = {},
  year    = {},
  pages   = {127--162},
}

@article{BML,
  author  = {E. Boros and MC Golumbic and VE Levit},
  title   = {On the number of vertices belonging to all maximum stable sets of a graph},
  journal = {Discrete Applied Mathematics},
  volume  = {124},
  year    = {2002},
  pages   = {17--25},
}

@article{PQ,
  author  = {JC Picard and M Queyranne},
  title   = {On the structure of all minimum cuts in a network and applications},
  journal = {Mathematical Programming Study},
  volume  = {13},
  year    = {1980},
  pages   = {8--16},
}

@article{HV,
  author  = {MA Henning and JH van Vuuren},
  title   = {Graph and Network Theory. An Applied Approach using Mathematica},
  journal = {Springer Optimization and its Applications},
  volume  = {193},
  year    = {2022},
  pages   = { 766 pages},
  }

\end{document}